\documentclass[a4paper]{article}
\usepackage{natbib}

  \usepackage{graphics,graphicx}
   \usepackage{amsmath}
   \usepackage{amsfonts}
   \usepackage{float,txfonts,amssymb,bm}
 \usepackage{colortbl}

\begin{document}

\title{An Intuitive Curve-Fit Approach to Probability-Preserving Prediction of Extremes.}
\author{Allan McRobie \\
Cambridge University Engineering Department\\
Trumpington St, Cambridge, CB2 1PZ, UK \\
fam20@cam.ac.uk}

\maketitle

\begin{abstract}
A method is described for predicting extremes values beyond the span of historical data.
The method - based on extending a curve fitted
to a location- and scale-invariant variation of the double-logarithmic QQ-plot - is simple and intuitive, yet it preserves probability to a good approximation. The procedure is developed on the Generalised Pareto Distribution (GPD),
but is applicable to the upper order statistics of a wide class of distributions.
\end{abstract}

\section{Introduction}
We present a method for extrapolating to extreme values outside the span of historical data.
The algorithm is accompanied by a visual representation -  a location- and scale-invariant variation of the double logarithmic QQ-plot - that accords with intuition.
A frequentist approach is taken, but the resulting predictor resembles one that might result from
a Bayesian predictive distribution with uninformative priors.

The approach begins  - via the curve-fit -  with an estimate $\hat{\xi}$ for the tail index $\xi_{DA}$, this being the tail (or shape) parameter of the Generalised Pareto Distribution (GPD) within whose domain of attraction the distribution lies. Predictions are then made from a GPD with an augmented tail parameter $\xi_p = \hat{\xi}+d\xi$, where the increment $d\xi$ depends on the estimate $\hat{\xi}$ and the desired recurrence level $T_{des}$ of the prediction. The delivered recurrence level $T_{del}$ is shown to be close to that desired for data drawn from a wide class of general distributions.

An extreme value predictor with good probability preservation was described in \cite{McRobiePRED}. That predictor was designed to return the highest point in the historical data sample of length $N$ as the level $T= N+1$ prediction, and more extreme predictions then formed a continuous curve emanating from the highest data point. In the paper here, the requirement for the predictor to pass through the highest historical data point is relaxed.

\section{The Basic Construction}

Whilst recognising that extrapolation to future events is a prediction problem, we nevertheless approach the problem
via estimation. We begin in the three parameter $ \lbrace \mu, \sigma, \xi \rbrace $ family of Generalized Pareto Distributions.
Estimation of the location and scale parameters $\mu$ and $\sigma$ is obviated by normalising the upper order statistics via ratios of data spacings. The shape parameter $\xi$ is estimated via a simple curve fit to a location- and scale-invariant variation of the double-logarithmic QQ-plot. Although there are tail estimators with marginally smaller mean square error, the curve fit estimator has explanatory benefits.

The GPD has tail distribution function $G = 1-F$ (where $F$ is the cdf)
\begin{equation}
G(x) = 1- F(x) = \left[ 1+ \xi \frac{x-\mu}{\sigma}\right]^{- 1/\xi}
\end{equation}

Consider $N$ data points ${x_1, \ldots x_N}$ sampled from a three parameter $ \lbrace \mu, \sigma, \xi \rbrace $  GPD ordered such that $x_1$ is the sample maximum, and consider just the first $k$ upper order statistics.

The data can be made location- and scale-invariant using the $j$th and $k$th order statistics. That is, for each data point $x_i$ we define
\begin{equation}
u_i = \frac{x_i - x_j}{x_j - x_k} \ \ \  \text{    with } j < k \leq N
\label{scaleddata}
\end{equation}

\begin{figure} [t!] \centering
  \hspace*{15mm}
 \includegraphics[width=160mm, keepaspectratio]{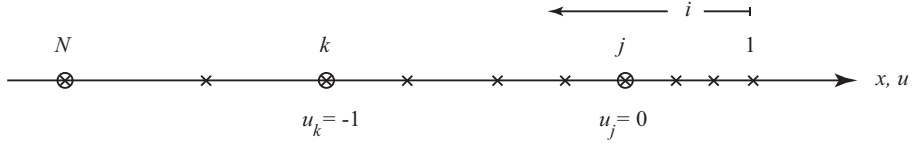}
\caption{Ordering and scaling the data}
\label{fig1}
\end{figure}

We now create the curves of the basic construction. Inverting the tail distribution function $G(x)$ of the GPD gives
\begin{equation}
\tilde{x}(G) = \mu + \frac{\sigma}{\xi}\left( G^{-\xi} -1 \right)
\end{equation}
and for given $j$ and $k$ we define
\begin{equation}
\tilde{u}_i = \frac{\tilde{x}(G_i) - \tilde{x}(G_j)}{\tilde{x}(G_j) - \tilde{x}(G_k)}
\end{equation}
We then choose abscissa $X$ and ordinate $Y$ as
\begin{eqnarray}
X & = & \log(1 + \tilde{u}_i) \\
Y & = & -\log( 1 - V_i) \ \ \text{  with   } \ \   V_i =   \frac{\log(G_i / G_j)}{\log(G_k/G_j)}
\end{eqnarray}

\begin{figure} [t!] \centering
 \includegraphics[height=105mm, keepaspectratio]{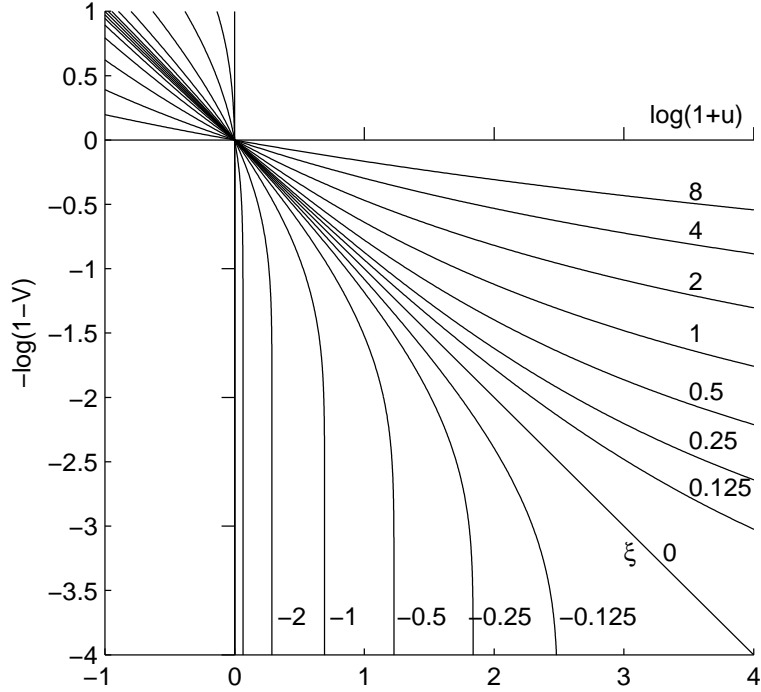}\\
\caption{The basic construction. Loosely speaking, the $X$ axis is the logarithm of the normalised data and the $Y$ axis is the double logarithm of the exceedance probability. }
\label{basic}
\end{figure}

We now make the approximation $G_i  \approx i/(N+1)$, and choose to have $k$ even and $j = k/2$. The $X$ coordinate thus depends on \begin{equation}
\tilde{u}_i   =  \frac{G_i^{-\xi} - G_j^{-\xi}   }{G_j^{-\xi} - G_k^{-\xi}   }  = \frac{ g_i^{\xi} -1  }{ 1 - \alpha^{\xi} } \label{utilde}
\end{equation}
with $g_i \equiv G_j/G_i = k/(2i)$ and $\alpha \equiv G_j/G_k = 1/2$, and the Y-coordinate becomes
\begin{equation}
Y = -\log\left[-\log(i/k)\right] + \log \log 2  \label{Yloglog}
\end{equation}
Noting that both $X$ and $Y$ are independent of $N$, we can plot, for each $\xi$, a curve parameterised by $r = i/k$ for all $0 < r \leq 1$. This gives the basic construction, shown in Figure~\ref{basic}. Loosely speaking, the $X$ axis is the logarithm of the normalised data and the $Y$ axis is the double logarithm of the exceedance probability.

The construction has a number of features:
\begin{itemize}
\item in the lower right quadrant, the GPD with $\xi$ = 0 (the exponential distribution) plots to the falling diagonal. The heavy-tailed GPDs ($\xi$ positive)
 plot above this, and the truncated tail GPDs ($\xi $ negative) plot below it.
 \item the point corresponding to $i = j = k/2$ plots to the origin, with larger data points plotting into the lower right quadrant.
 \item for the GPD, the construction is independent of $N$.
\end{itemize}

Data drawn from a GPD can now be plotted on this diagram. The data defines $X$ coordinates, via the statistic $u$ (Eqn.~\ref{scaleddata}), and the ordinates $Y$ are defined solely by the order index ratio $r = i/k$ (via Eqn.~\ref{Yloglog}).  Guided more by intuition than theorems, the general idea is that data from a GPD with tail parameter $\xi$ is in some sense likely to plot near to the underlying curve defined by $\xi$, such that an estimate of $\xi$ can be obtained by, say, a simple least-squares curve fit.

\begin{figure} [h!] \centering
   \includegraphics[height=55mm, keepaspectratio]{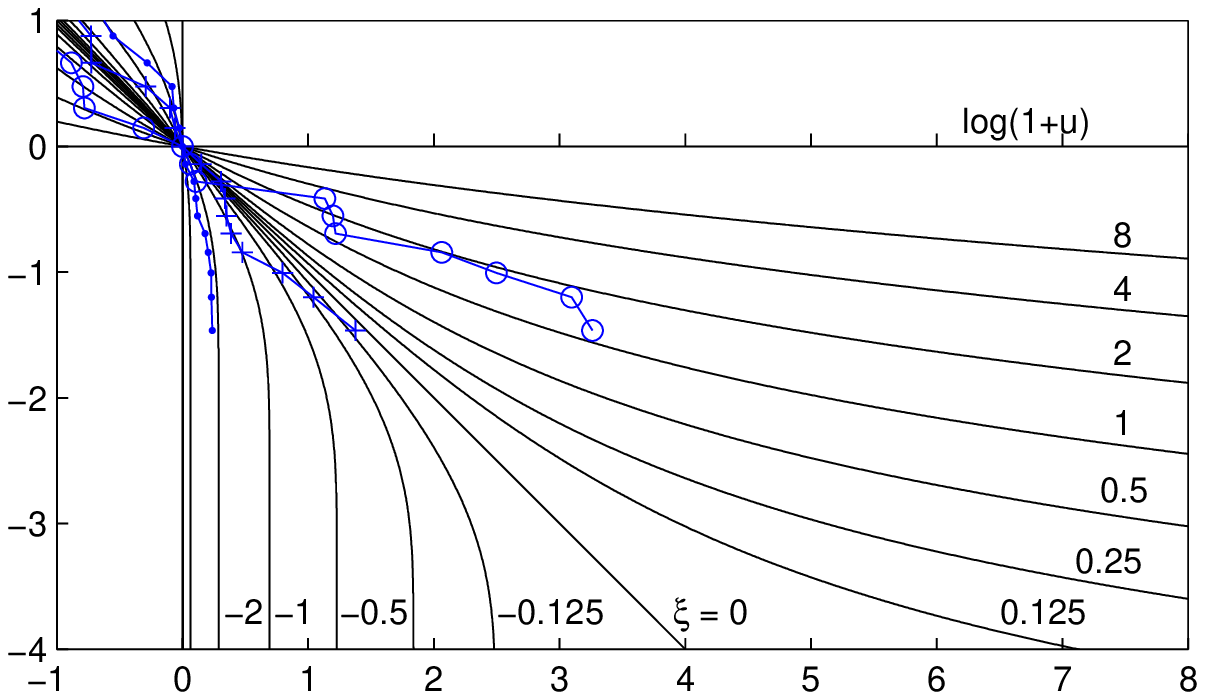}
  \includegraphics[height=55mm, keepaspectratio]{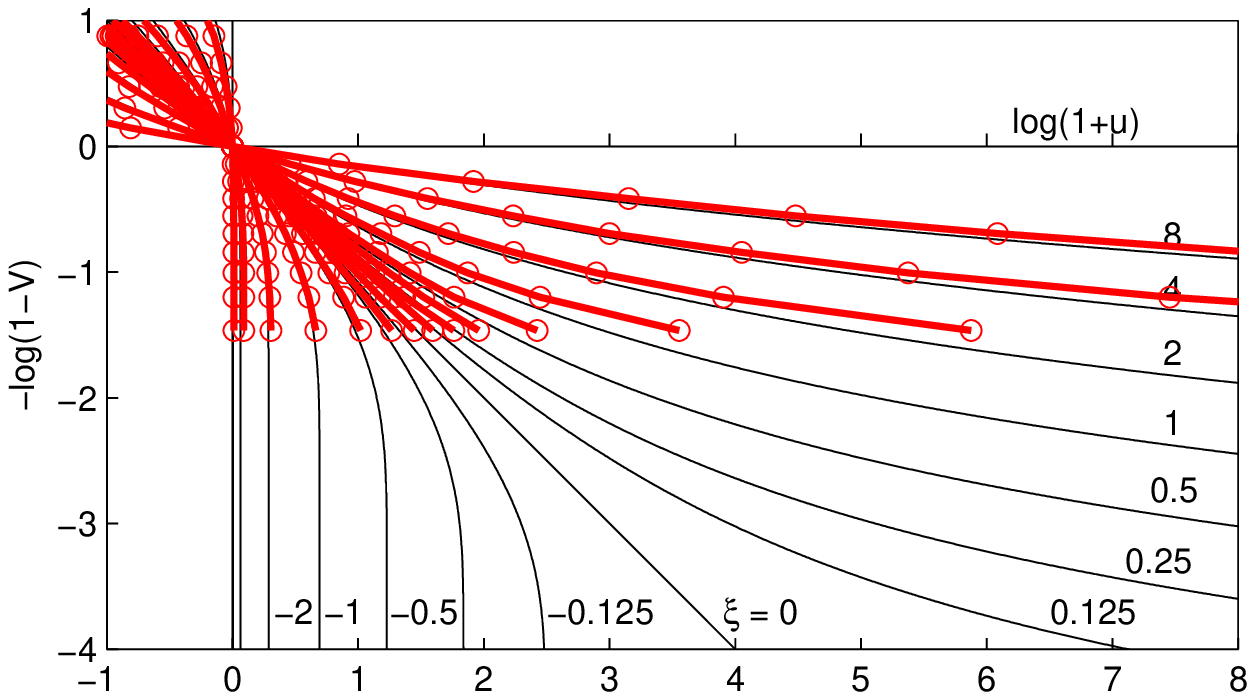}\\
    \includegraphics[height=55mm, keepaspectratio]{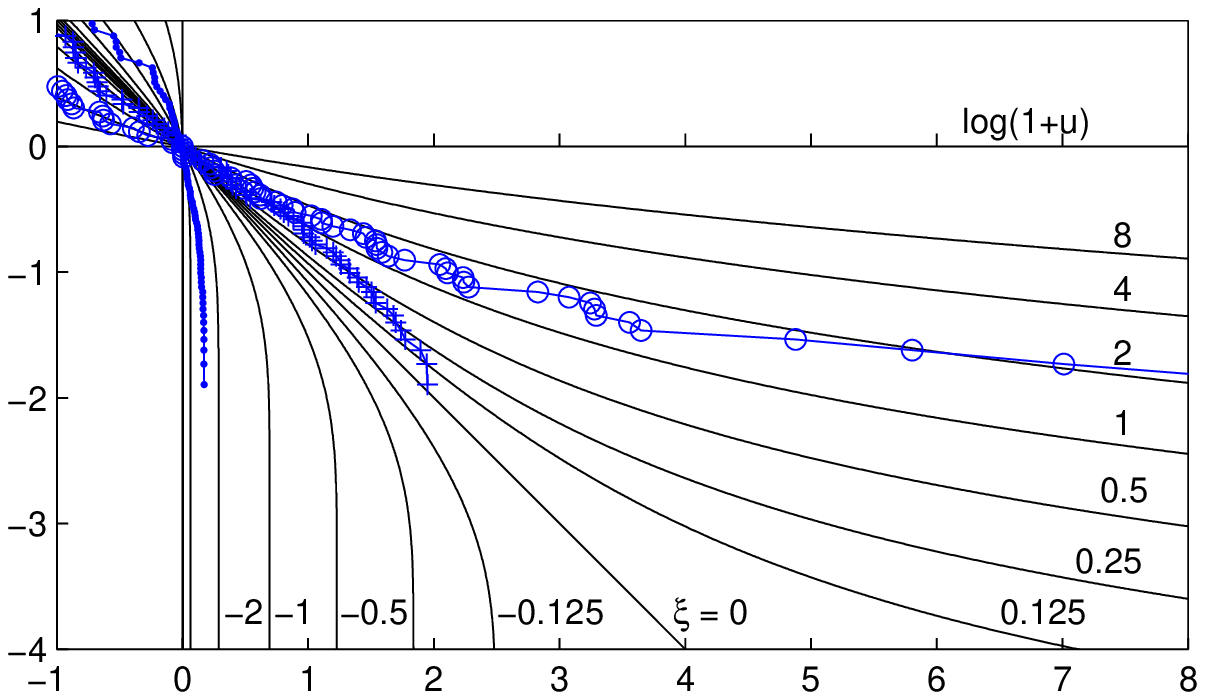}
  \includegraphics[height=55mm, keepaspectratio]{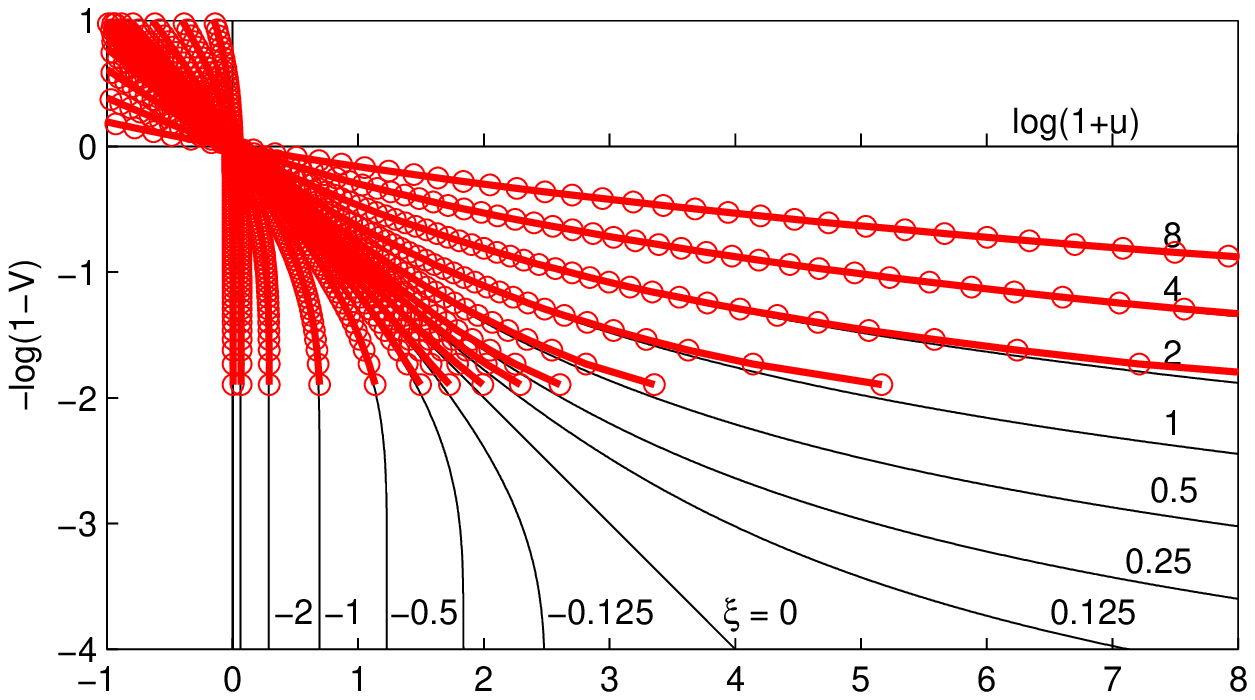}
 \caption{Plotting data onto the basic construction. The left-hand diagram shows specific GPD samples at $\xi = -2,0,2$ (labelled $\cdot$ + o respectively) and the right-hand diagram shows the average results for 10000 samples at each of the $\xi$ labelled. The upper and lower diagrams correspond to $k= N = 20$ and $k=N=200$ respectively.}
\label{basics}
\end{figure}

Figure~\ref{basics} shows individual samples (left) and sample averages (right). The individual samples lie in the approximate vicinity of their respective curves. The sample averages (that is, the averages of $\log(1+u_i)$ over 1000 samples) follow their corresponding curves reasonably closely, although with a tendency to follow a curve at a slightly higher $\xi$ than that from which they were sampled. This tendency appears to reduce as the sample size $N$ is increased ($N=20$ in the upper figures, and $N = 200$ in the lower figures).

The intuitive visual representation afforded by such figures can be helpful in indicating how tail estimators applied to data drawn from more general non-GPD distributions can often lead to consistently erroneous estimates. We give an example in the Appendix, noting here only that - as we are generally trying to estimate the parameter that represents the shape of the tail - the ability of the $\log \log$ versus $\log$ construction to make visible the shape of the tail can be instructive.

Alternative plotting positions can be obtained by using $G_i = (i-0.5)/N$, as is commonly done in QQ-plots  (see for example \cite{embrechts}, p292). For a traditional empirical distribution function, this would correspond to the midpoint of the step at each data point. With these new plotting positions, the previous equations still apply, but with the revised definitions for the $G$'s included throughout.

The definition of the plotting position has a significant effect on where data is plotted. It also affects the background analytical curves. However, whilst points on a curve may move significantly, the curve as a whole moves little. Since the graphic is for illustration only, we may therefore use any $k$ to draw the background curves, but adopt the appropriate exact values when doing computations.

\begin{figure} [h!] \centering
   \includegraphics[height=55mm, keepaspectratio]{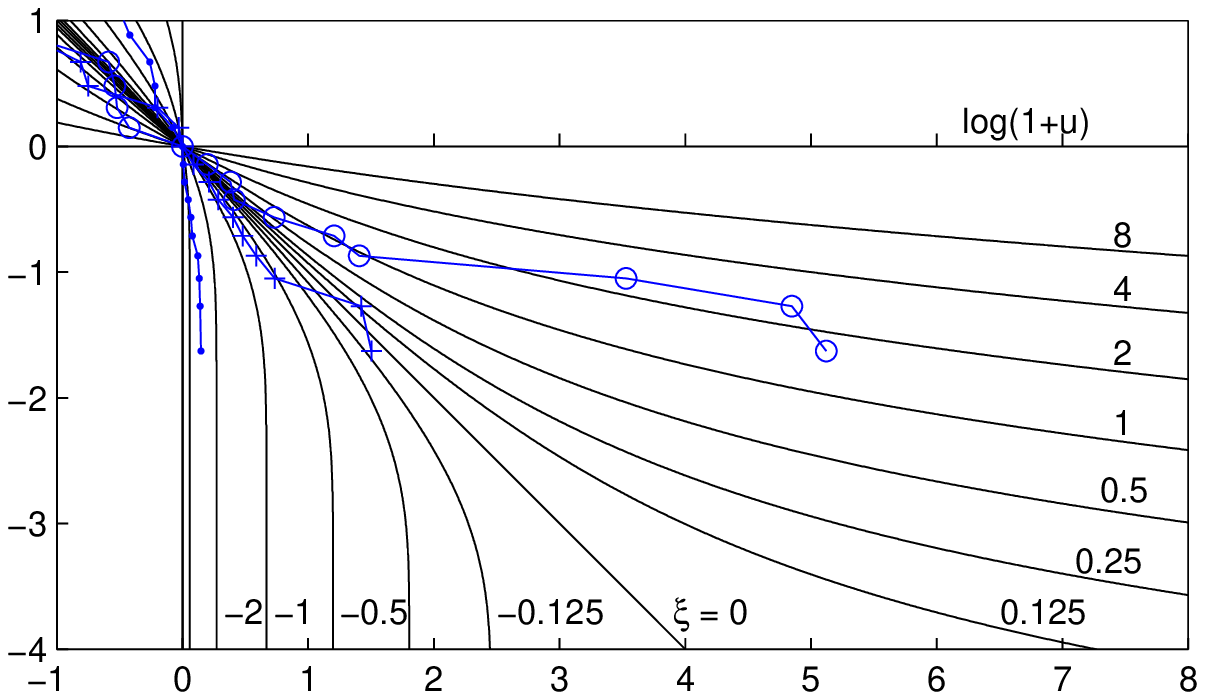}
  \includegraphics[height=55mm, keepaspectratio]{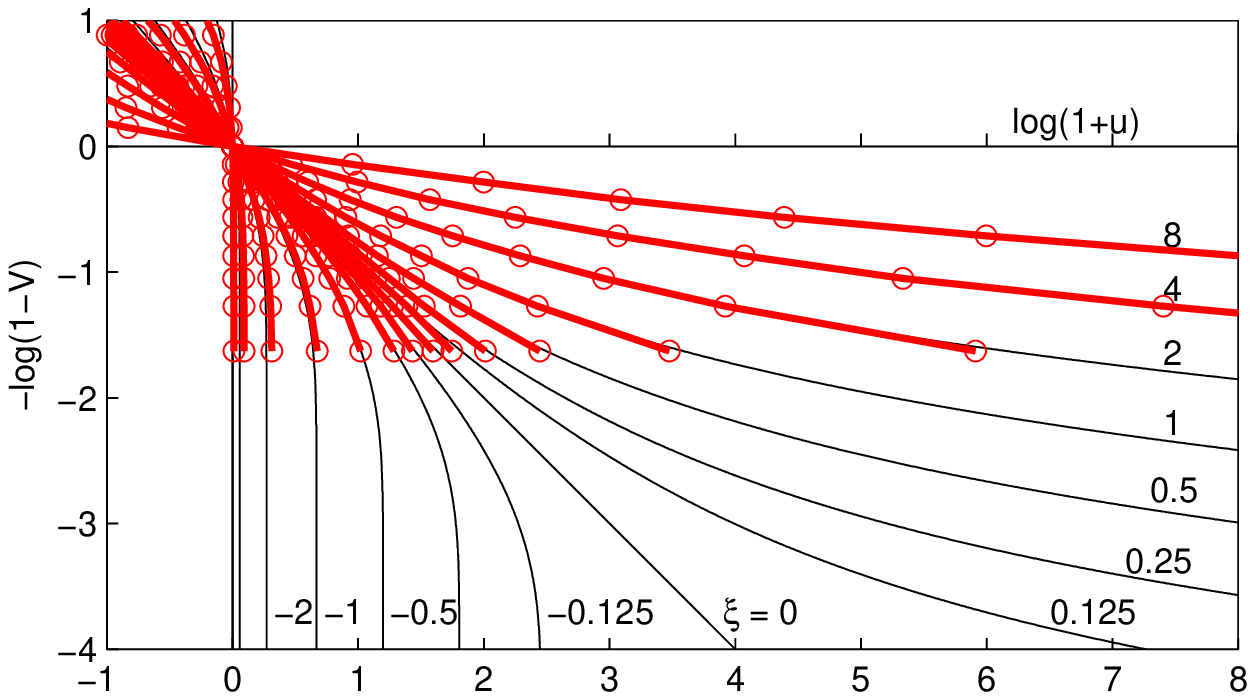}
 \caption{The $(k- 0.5)$ construction. Specific GPD samples ($N=k=20)$ at
 $\xi = -2,0,2$ are shown left, and the averages for 10000 samples at each of a range of $\xi$ are shown right.}
\label{khalf}
\end{figure}

We shall call such diagrams the $(k-0.5)$ construction, and examples are shown in Figure~\ref{khalf}. There, GPD data has been added and the averages of sampled data now plot much closer to the (re-drawn) curves of the underlying analytical approximations. We shall thus use the $(k-0.5)$ construction for the estimation phase, and apply a least-squares curve-fit to estimate $\xi$.

(Unless stated otherwise, all further GPD samples and analysis will use $N=k=20$.)

\subsection{Curve Fitting}

The ordinates are fixed by the plotting position definitions.
We therefore consider the horizontal error between data and curve, defining $\epsilon \equiv f(u_i) - f(\tilde{u}(G_i))$ for some function $f$, and  we choose $f(u) = X =  \log(1+u)$ to accord with the construction.

The sum of the squared errors over the first nine data points (the tenth data point having no error by this construction) is
\begin{equation}
 \sum_i \epsilon_i^2 = \sum_{i=1}^{9} \left( f(u_i) - f(\tilde{u}_i) \right)^2
\end{equation}
where $u_i$ is the scaled data and $\tilde{u}_i \equiv \tilde{u}(G_i)$ is its analytical approximation, with $G_i = (i-0.5)/N$.

Taking the derivative with respect to the tail parameter gives
\begin{equation}
\frac{d}{d\xi} \sum_i \epsilon_i^2 = 2 \sum_i \left( f(u_i) - f(\tilde{u}_i) \right) \ f'(\tilde{u}_i) \ \frac{d\tilde{u}_i}{d\xi}
\end{equation}

The analytic approximation for the $i$th order statistic is
\begin{equation}
\tilde{u}_i = \frac{g_i^\xi -1}{1- \alpha^\xi}
\text{     with derivative   }
\frac{d\tilde{u}_i}{d\xi} = \left( \frac{1}{1-\alpha^\xi}\right)\left[ g_i^\xi \log g_i + \tilde{u}_i \alpha^{\xi} \log \alpha  \right]
\end{equation}
Since we choose $f(u) = \log(1+u)$  we have $f' = 1/(1+u)$, such that error minimisation requires
\begin{equation}
\frac{d}{d\xi} \sum_i \epsilon_i^2 = 0 = \sum_i  \frac{1}{1+\tilde{u}_i} \ \frac{d\tilde{u}_i}{d\xi} \ \log \frac{1+u_i}{1+\tilde{u}_i}
\label{curvefiteqn}
\end{equation}
This can be readily solved numerically for the estimate $\hat{\xi}$ (using the Matlab $fzero$ function, for example)

\begin{figure} [ht!] \centering
   \includegraphics[height=50mm, keepaspectratio]{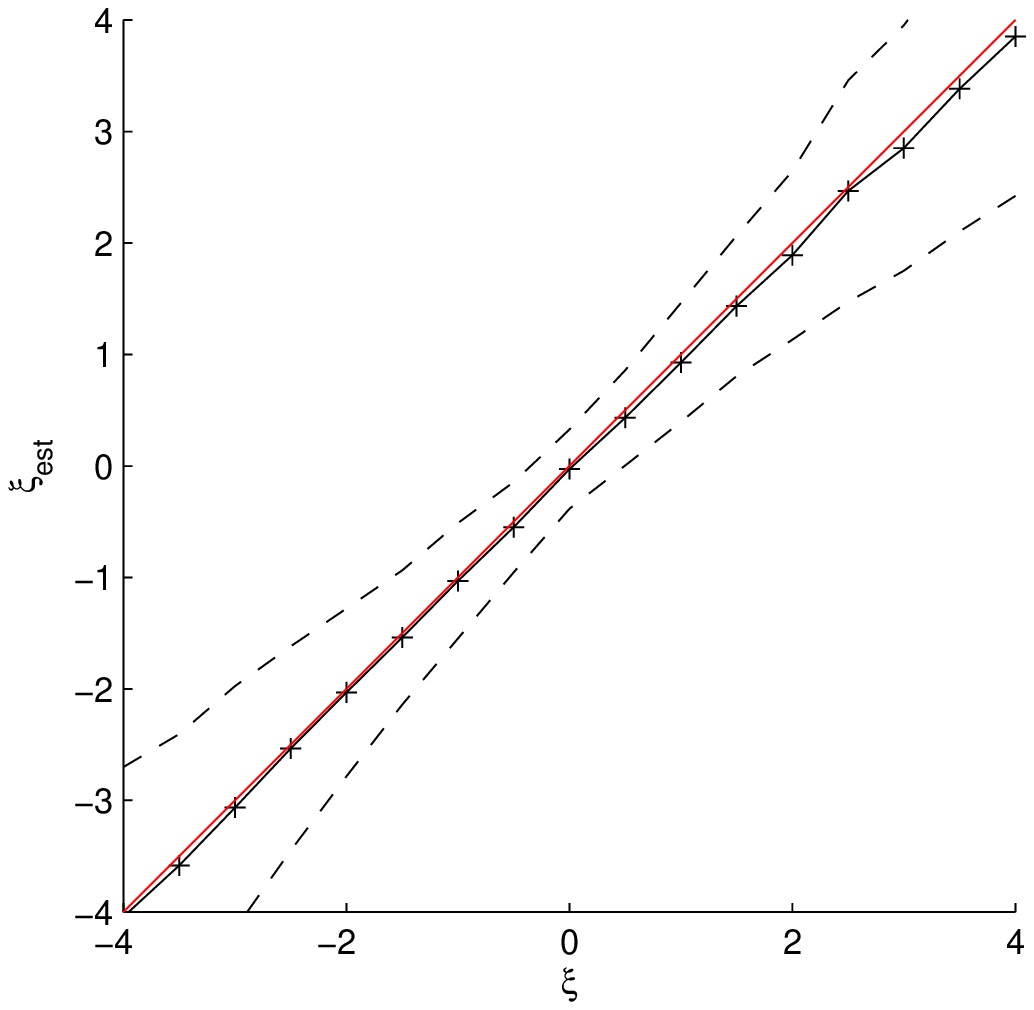}
  \includegraphics[height=50mm, keepaspectratio]{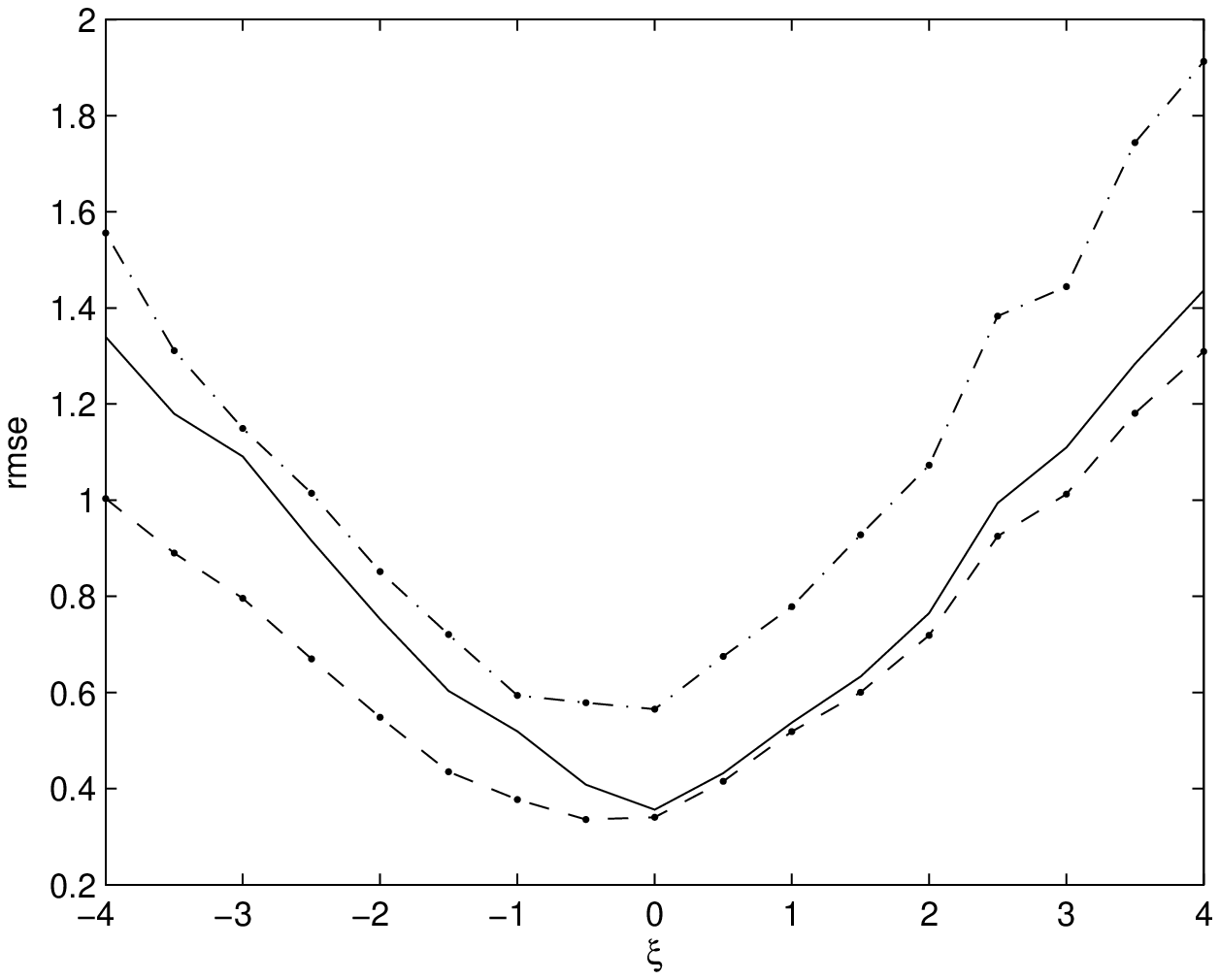}\\
 \caption{The left diagram shows the performance (mean, and mean $\pm$ std) of the curve-fit estimator for 1000 samples at each of a wide range of $\xi$. The right diagram shows the root mean square error of the curve-fit estimator (solid), nestling between the rmse performance of the ``linearly-rising'' elemental estimator that uses  all 20 data points (lower curve) and the one that uses just the first ten (upper curve).}
\label{fig3}
\end{figure}

The performance of the resulting estimator is shown in Figure~\ref{fig3} for samples drawn from GPDs over a range of $\xi$. The left figure shows that the estimator has some variance about a very small bias, and the right figure compares the root mean square error to that of an alternative estimator - the ``linearly rising'' combination of elemental estimators
described in \cite{McRobieGPD}. The curve-fit estimator is based on the first ten upper order statistics, together with knowledge of the twentieth. Two versions of the elemental estimator are shown in the right-hand figure. The upper curve uses knowledge of only the first ten upper order statistics, whereas the lower (more accurate) curve uses knowledge of all twenty. The curve-fit estimator nestles between the two.

Although there exist tail estimators with smaller rmse (e.g. \cite{segers}) we proceed using the curve-fit for a number of reasons.  First, despite its simplicity, its performance is only marginally worse than the best available alternatives, and in many cases it is superior to other commonly-adopted procedures. For example, the curve fit works for all $\xi$, not just $\xi$ positive (as is the case for Hill's estimator), it does not suffer from the numerical difficulties encountered by Maximum Likelihood estimation (see for example \cite{castillo, McRobieGEV}) nor does it have the convergence issues associated with MCMC Bayesian approaches applied to extremes. Secondly, only the most basic description of the curve-fit is described here and there is considerable scope for further optimisation later.
Finally - and most importantly -  curve-fitting brings considerable intuitive insight to the problem. This will become evident as we proceed from estimation to prediction, and especially in the case of non-GPD data.

\section{Prediction in the GPD case}
\begin{figure} [ht!] \centering
   \includegraphics[height=50mm, keepaspectratio]{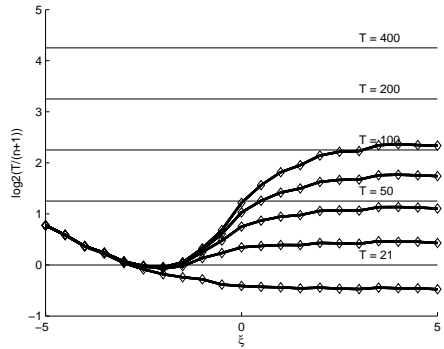}
 \caption{Performance of the naive predictor that predicts out of sample extremes using the GPD with tail parameter set to the value $\hat{\xi}$ that was estimated from a sample of size $N=20$ drawn from a GPD with parameter $\xi$. The five curves show the return levels delivered at various $\xi$, and the horizontal lines are the target return levels. It can be seen that this approach has a tendency to underpredict. (The slight overprediction for $T= 21$ at the lower left is thought to be an inconsequential result of the plotting position redefinition). }
\label{fig4}
\end{figure}

Having obtained an estimate $\hat{\xi} $ of the tail parameter $\xi$, a naive approach to prediction beyond the span of the data would be to extrapolate the GPD with parameter $\hat{\xi}$ out to large extreme values. However, as is well known, this is optimistic and is liable to considerably underestimate possible future extremes. This is illustrated in Figure~\ref{fig4}. Here, for each of 10000 samples of size $N = 20$ drawn from a pure GPD at each of a variety of values of $\xi$,  the $k=20$ curve-fit estimate $\hat{\xi}$ is first calculated. Predictions are then made at the recurrence levels $T = 21$, 50, 100, 200 and 400 using
\begin{equation}
x_T = x_j+(x_j-x_k)u_T
\label{predeqn1}
\end{equation}
 with the normalised prediction given by
\begin{equation}
u_T  = \frac{g_T^{\hat{\xi}} -1}{ 1 - \alpha^{\hat{\xi}}}
\label{predeqn2}
\end{equation}
where $g_T = G_j/G_T$. Here the recurrence probability of the extreme event at level $T$ is set to $G_T  = 1/T$ and for this prediction phase we return to the plotting position approximations $G_j = j/(N+1)$ and $G_k = k/(N+1)$ in line with the original basic construction.

To test each prediction, a further 10000 data points are sampled from the same distribution, and the number $N_{exc}$ that exceed the prediction is counted. For each sample and its associated prediction, the exceedance probability is thus estimated as  $G_{del} = N_{exc}/10000$. The procedure is repeated for 10000 samples and the delivered recurrence level is estimated as $T_{del} = 1/\text{mean}(G_{del})$. (Note: because the form of the supposedly ``unknown'' distribution from which data is sampled is actually known in simulations, $G_{del}$ could be evaluated analytically for each prediction. However, the numerical method is more readily adaptable and leads to results which are substantially the same).

Figure~\ref{fig4} shows the delivered ``return period'' to be generally substantially lower than the desired return level (except for $T=21$ at the far left - and this latter is thought to be a minor and irrelevant consequence of using the $(i-0.5)/N$ plotting positions for estimation). The wider tendency of naive extrapolation to underpredict can be readily appreciated by consideration of samples drawn from a GPD with say $\xi = 0$. For around half the samples, the tail parameter will be estimated as negative, corresponding to a bounded tail distribution, and extrapolating to rare events in the bounded tails leads to predictions which are readily overcome by further samples drawn from the true $\xi = 0$ distribution.

To compensate for such underprediction, the {\bf central proposal of this paper} is to predict instead from GPDs with tail parameters higher than estimated. That is, having obtained the estimate $\hat{\xi}$, prediction is then accomplished using
$\xi_{p} = \hat{\xi} + d\xi$, where $d\xi$ depends on both the estimate $\hat{\xi}$ and the desired return level $T_{des}$.
In graphical terms, having done a GPD curve-fit to the data on the $(k-0.5)$ construction (as per Figure~\ref{khalf}), prediction to more extreme return levels proceeds along a curve that increases to values progressively to the right of that originally estimated. The curve so constructed corresponds, in some sense to the existence of a ``predictive distribution''. Although predictive distributions are well defined in the Bayesian framework, they appear to be little used and somewhat ill-defined in frequentist settings.

\subsection{The increment $d\xi$}
Trial and error numerical experiments led to the construction of a set of curves for the increment $d\xi$. That is, at some given desired return level $T_{des}$, a function $d\xi$ was defined over the complete range of $-\infty < \hat{\xi} < \infty$. GPD samples were generated at each $\xi$ over a wide range ($-5 \leq \xi \leq 5$), and for each sample, a prediction was made using the GPD with $\xi_p = \hat{\xi}+d\xi$ as per Eqns.~\ref{predeqn1} and \ref{predeqn2} (but with the latter using $\xi_p$). The number of times the prediction was exceeded by further samples drawn from the GPD at the original $\xi$ was recorded, allowing the delivered recurrence level $T_{del}$ to be estimated. If $T_{del}$ fell below $T_{des}$ over particular ranges of $\xi$ the function $d\xi$ was increased smoothly and locally in that general area. Typically the process began by calibrating an exponential function to get good probability matching at extreme negative and positive $\xi$ and then augmenting it with a set of ad hoc bump functions - broad Gaussians - in the regions of moderate $\xi$. Regions of $\xi$ moderate and negative typically required most augmentation.
The bump functions were adjusted manually until $T_{del} \approx T_{des}$ for all $\xi$ tested $(-5\leq \xi \leq 5)$. The process was then repeated to determine the function $d\xi$ at another recurrence level $T$.

\begin{figure} [ht!] \centering
   \includegraphics[height=120mm, keepaspectratio]{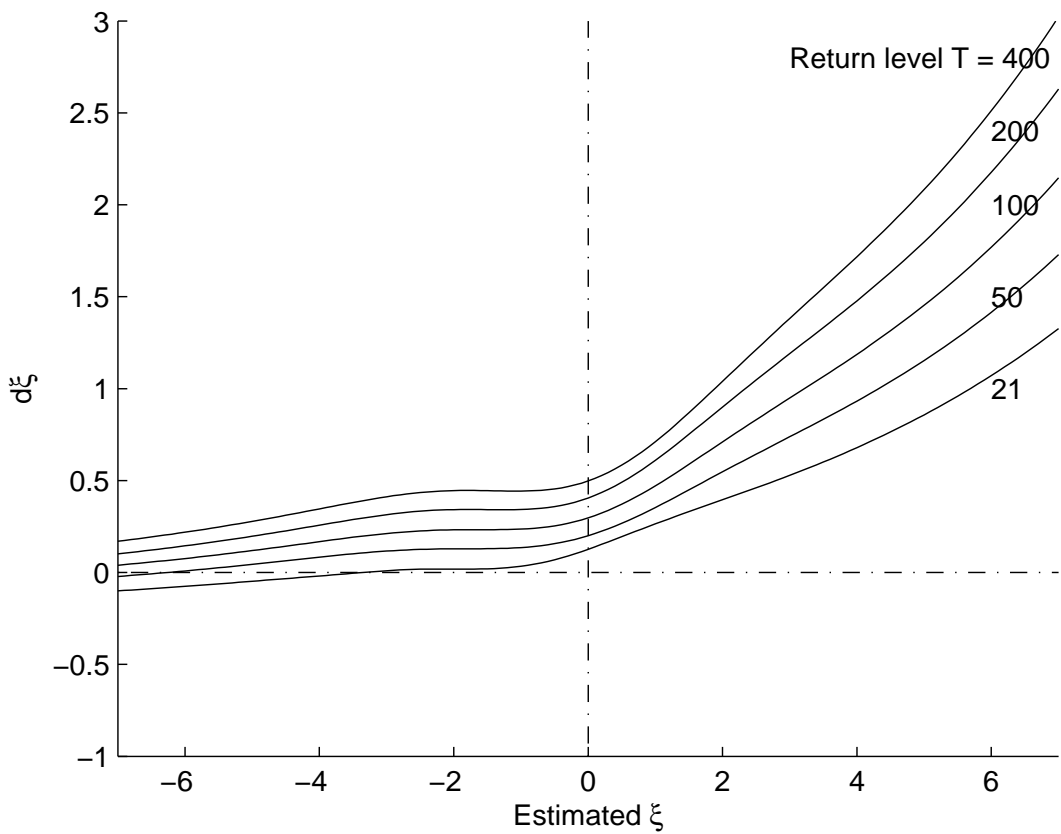}
 \caption{The increment $d\xi$ to be added to the estimate $\hat{\xi}$ to obtain predictions with good probability preservation.}
\label{fig5}
\end{figure}

The resulting functions $d\xi(\hat{\xi}, T_{des})$ are shown in Figure~\ref{fig5} for five recurrence levels $T_{des} = 21$, 50, 100, 200, 400,  these values being chosen for practical reasons. The functions have no neat analytic form.

\begin{figure} [ht!] \centering
  \includegraphics[height=50mm, keepaspectratio]{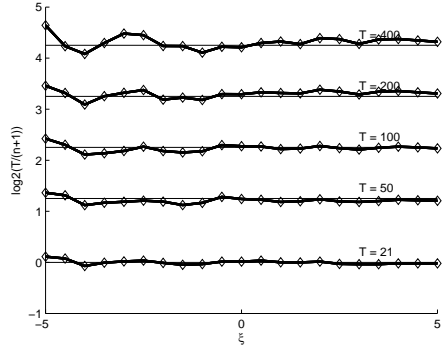}
 \caption{Performance of the predictor based on $\xi_{p} = \hat{\xi} + d\xi$. The delivered return levels are very close to target levels, even for extrapolations to return levels as large as $T= 400$ from $N=20$.}
\label{fig6}
\end{figure}

The delivered recurrence levels of the predictor are shown in  Figure~\ref{fig6}, where it can be seen that the performance is very close to target, even for desired return levels as high as $T_{des} = 400$, this being a significant extrapolation from $N=20$ data points. Without suggesting that extrapolating far outside the data is in any way advisable, if out of sample extrapolation is to be undertaken then a method that preserves probability (such as the predictor here does to a good approximation) would appear to be  more rational and less imprudent than naive extrapolation.

\begin{figure} [h!] \centering
   \includegraphics[height=130mm, keepaspectratio]{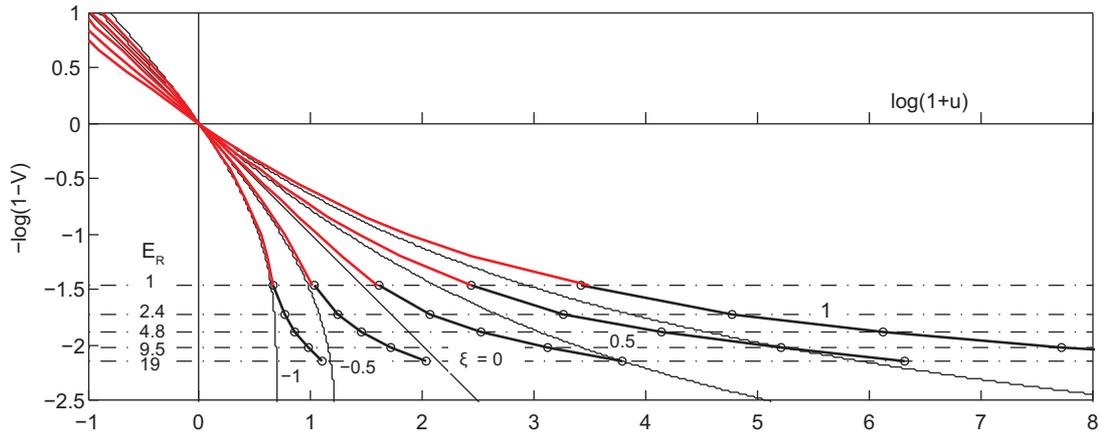}
 \caption{The basic extrapolation curves using $\xi_{p} = \hat{\xi} + d\xi$ (for $ \hat{\xi} = -1, -0.5, 0, 0.5 ,1$) are plotted in thick black at the lower part of the diagram. The diagram assumes the model is GPD with $N=k=20$, and extrapolations are shown to recurrence levels $T$ ranging from 21 to 400, corresponding to extrapolation ratios $E_R$ from 1 to 19. Above this is the part of the diagram usually used for estimation. Here the average values of 10000 samples of normalised data have also been plotted in red. Note that this diagram uses $G_i = i/(N+1)$ plotting positions throughout, in line with the basic construction.}
\label{fig7}
\end{figure}

The predictor can also be illustrated graphically on the basic construction. The standard analytical approximation curves can be plotted for a set of $\xi$. The region of the graph above the horizontal line corresponding to the return level $T = 21$ is the region used for estimation.  The corresponding extrapolation curves given by the $\xi_p = \hat{\xi}+d\xi(\hat{\xi},T)$ predictor can also be added below the horizontal $T= 21$  line. We thus obtain Figure~\ref{fig7}.

On this graph the averages of $\log(1+u_i)$ for data sampled from GPDs at the five values of $\xi$ have also been plotted in red, using $G_i = i /(N+1)$ plotting positions. The data averages thus now plot a little to the right of the approximate analytic curves (as was seen earlier in Figure~\ref{basics}). To reiterate, estimation uses $(i-0.5)/N$ and prediction uses $i/(N+1)$ plotting positions. The addition of the data here using the $i/(N+1)$ convention is thus for illustration purposes only. However, the way that the prediction curves below the $T = 21$ horizontal appear to match up so neatly with the data averages is remarkable, given that
the red curves correspond to data at a specific $\xi$, whilst the prediction curves are applied to data from GPDs of any $\xi$.
Even though, intuitively, some degree of matching may have been expected, there is as yet no explanation for this rather pleasing outcome that the prediction curves appear to be continuations of the data averages.

 The procedure for making a prediction from any data set is, of course, applied computationally, and the graphic of Figure~\ref{fig7} is for illustration only. Nevertheless, it does show how intuitive the method is.

 The way the prediction curves (thick black) veer to greater values than the naive extrapolation of the ``expected'' curves (thin black) is reminiscent of the behaviour of the Bayesian predictive distribution.  However the prediction curves here have been obtained  without reference to any prior distribution on the shape parameter. This is interesting, because there is as yet no known rational noninformative prior for the shape parameter in the three parameter GPD, and Bayesian extreme value theorists without access to any meaningful prior information may perhaps be inclined put broad normals centred loosely around $\xi = 0$. MCMC simulations also have computational issues for large $T$ predictions (associated with needing lengthy runs to populate the distant tail) whereas the procedure described here is essentially a simple function evaluation. In summary, that function evaluation is:

 \begin{itemize}
\item estimate $\xi \approx \hat{\xi}$ via a curve-fit to the scaled data (equation~\ref{curvefiteqn} as illustrated by Figure~\ref{khalf});
\item evaluate $d\xi(\hat{\xi}, T_{des} )$, the increment illustrated in Figure~\ref{fig5};
\item predict out-of-sample extremes using $\xi_p = \hat{\xi} + d\xi$ (as illustrated by Figure~\ref{fig7}).
 \end{itemize}

The algorithm has been described here in its most basic form, and there is considerable scope for further optimisation and extension.
For example, only the $N = k = 20$ case has been described, and extension to other values of $N$ and $k$ may be readily devised.
All normalisation has been with respect to the $(k/2)$th and $k$th order statistics, which again, can be generalised and optimised. Also, weighted least squares estimation could be used in the curve fit (the current weights being unity for $1 \leq i \leq 9$ and zero for $10 \leq i \leq 20$) and alternatives may lead to reduced estimation error. Similarly, least squares is applied to the function $\log(1+u_i)$ to accord with the basic construction, but other functions may reduce estimation error.

The ultimate objective, however, is less concerned with improving parameter estimation, and more with making out-of-sample extrapolations in the general domain-of-attraction case, which we now consider.

\section{Application to non-GPD data}
 Extreme Value Theory tells us that the upper order statistics of samples drawn from a non-GPD distribution may, with suitable scaling, in appropriate limits and under appropriate conditions, be approximated as a GPD.  The non-GPD distribution is then said to be in the domain of attraction of that GPD with tail parameter $\xi_{DA}$, and $\xi_{DA}$ is said to be the {\it tail index} of the non-GPD distribution. We now use this approximation to apply our GPD predictor to non-GPD data.

 The strategy here is fairly obvious. Given a sample of size $N \geq 20$ from a non-GPD distribution, select the $k=20$ upper order statistics and use least squares on the $(k-0.5)$ construction to make an estimate $\hat{\xi}$ of the tail index $\xi_{DA}$. Then make location- and scale-invariant predictions using the GPD with tail parameter $\xi_p = \hat{\xi}+ d\xi$. Finally, to demonstrate that the procedure works, the performance of the predictor will be tested by seeing how often the prediction is exceeded by further data drawn from the original non-GPD distribution.

In the pure GPD case we used $N=k=20$ and extrapolated to exceedance levels $T = 21$, 50, 100, 200, 400, corresponding to extrapolation ratios $E_R$ of $E_R = T/(N+1) = T/21 =$ 1, 2.4, 4.8, 9.5 and 19.   In the non-GPD case, we allow $N$ to be larger than $k=20$. Predictions however will still be at the same extrapolation ratios from the 20 upper order statistics, and the recurrence levels with respect to the full sample size $N$ will be $T = (N+1) E_R$.

\begin{figure}[ph!] \centering
\includegraphics[height = 55mm,keepaspectratio]{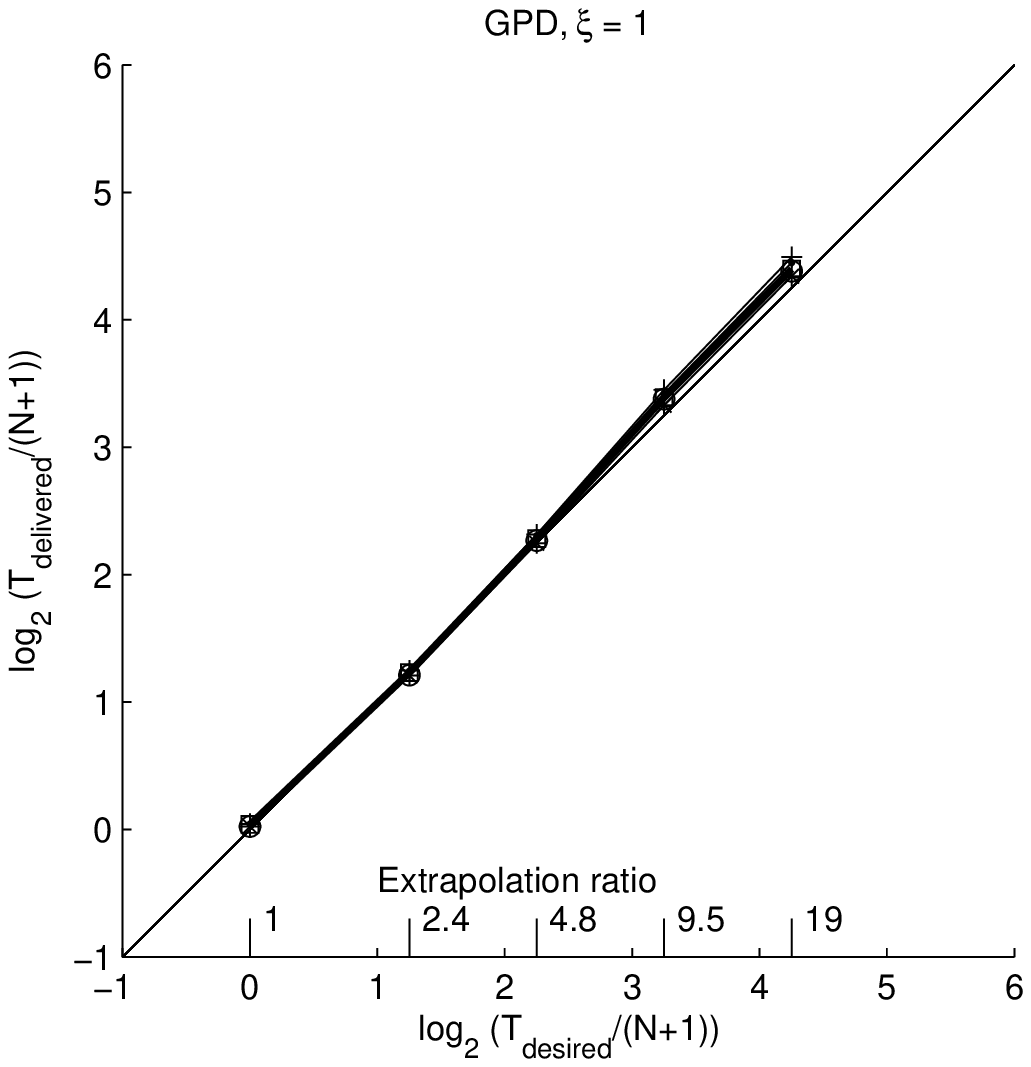}
\includegraphics[height = 55mm,keepaspectratio]{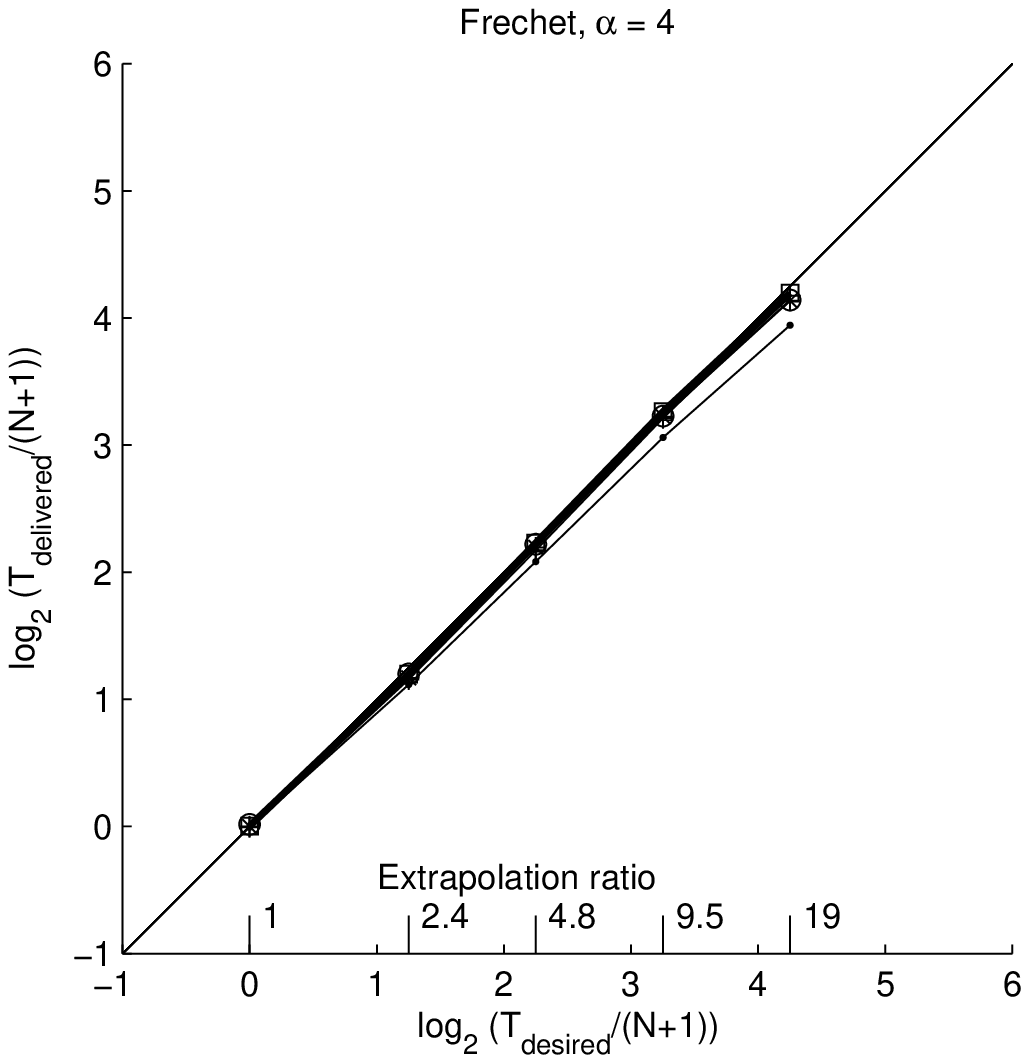} \\
\includegraphics[height = 55mm,keepaspectratio]{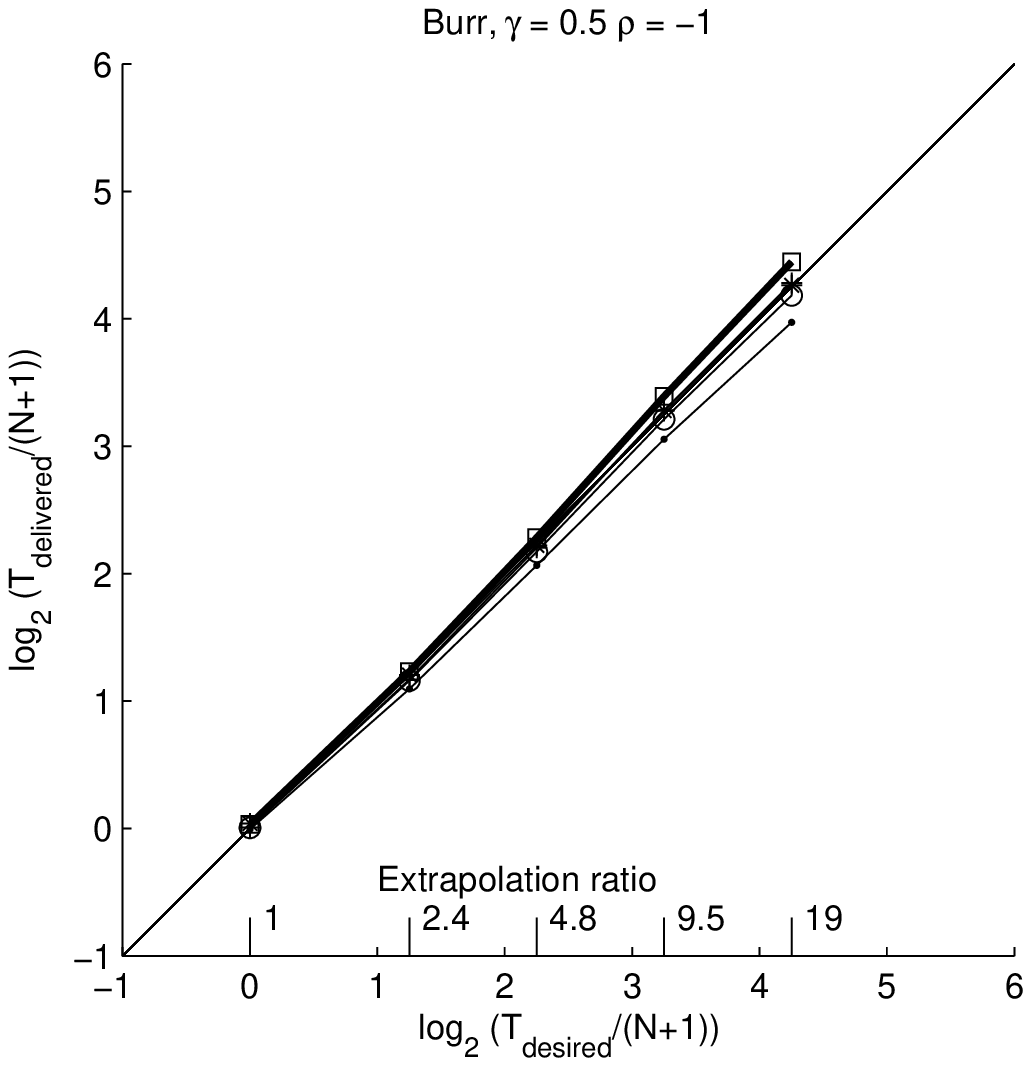}
\includegraphics[height = 55mm,keepaspectratio]{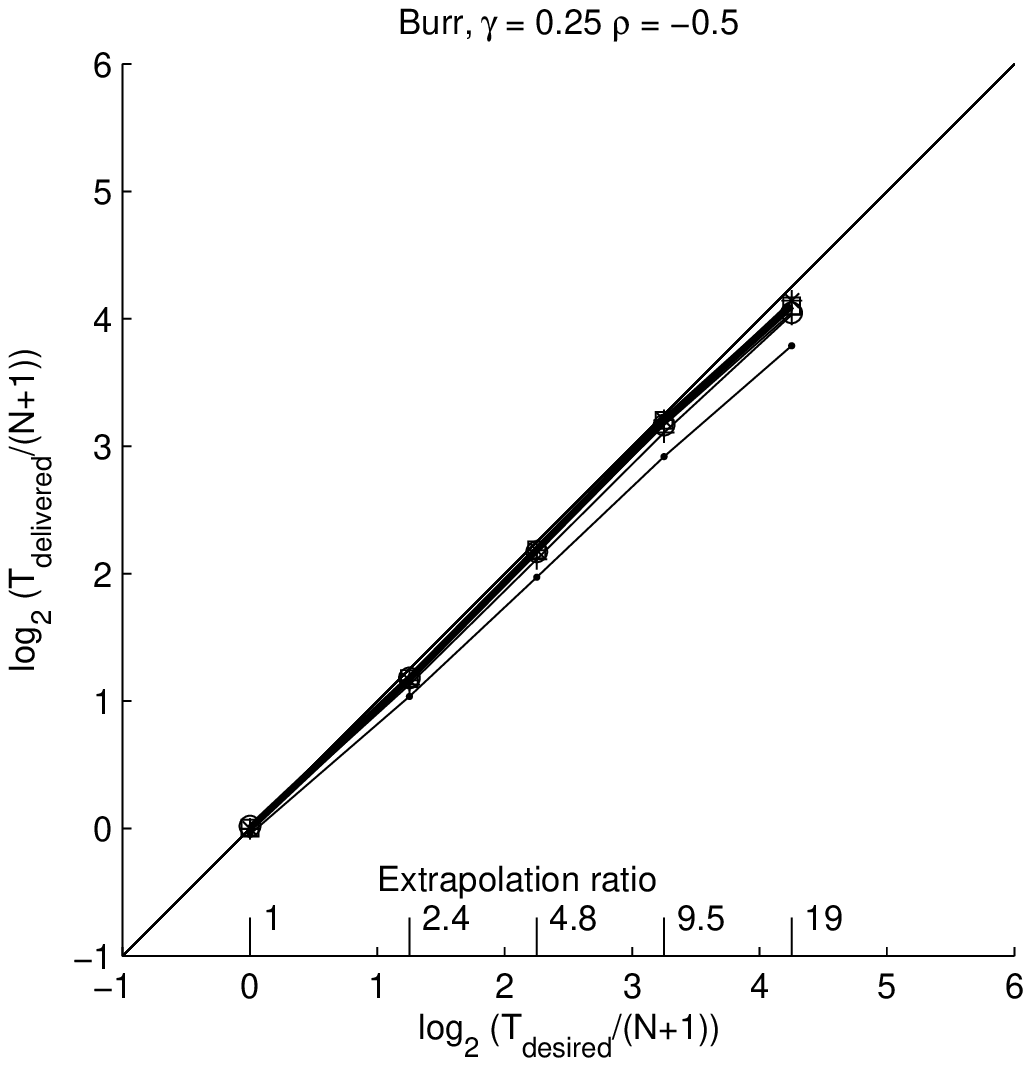} \\
\includegraphics[height = 55mm,keepaspectratio]{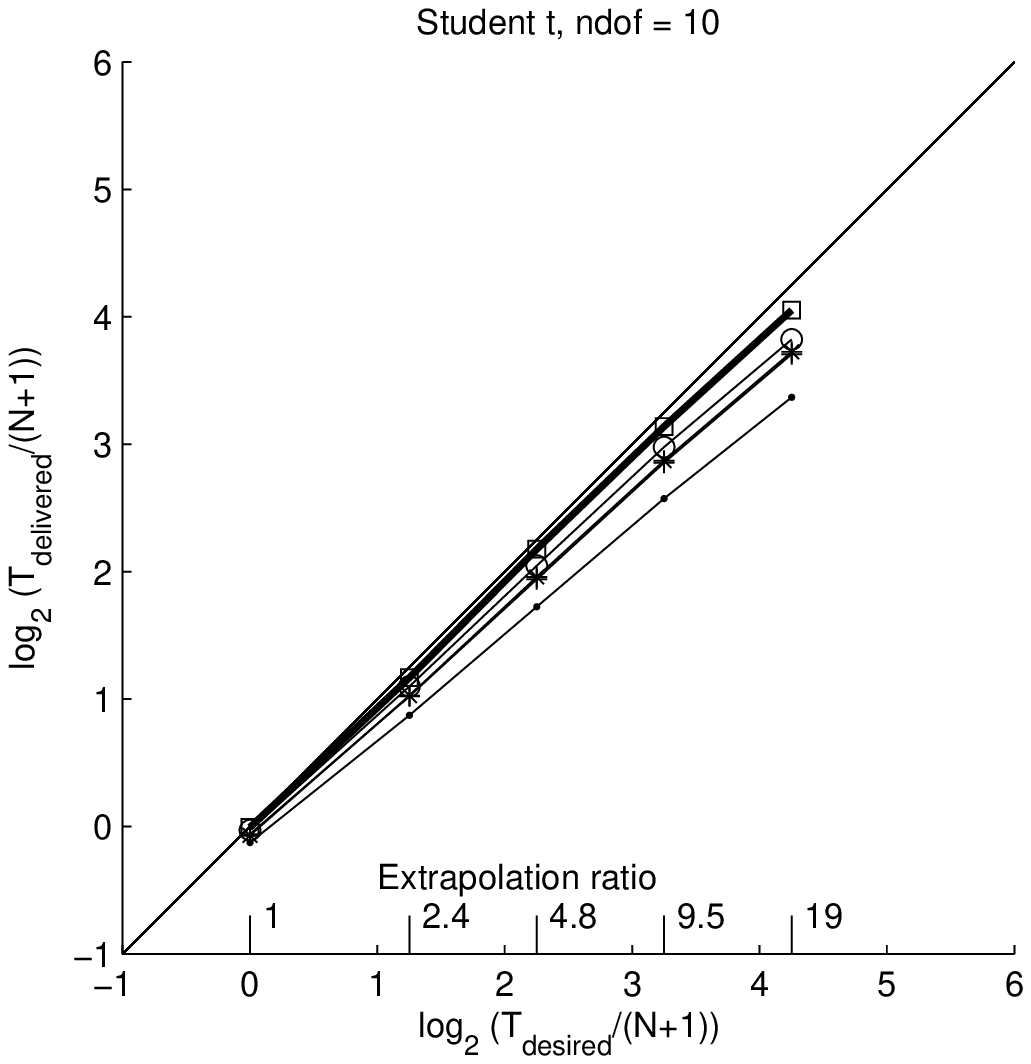}
\includegraphics[height = 55mm,keepaspectratio]{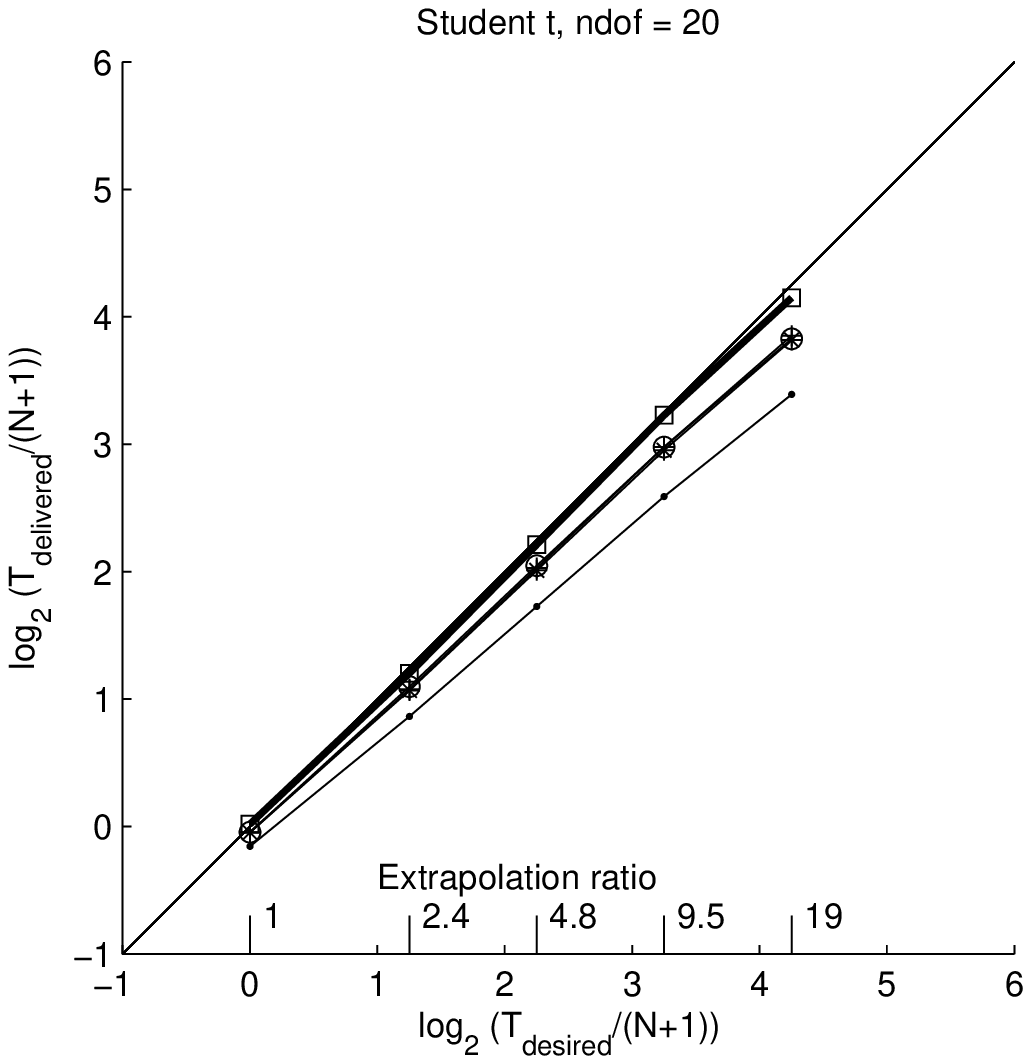}
  \caption{Probability performance of the predictor applied to samples drawn from GPD ($\xi = 1$), Frechet, Burr and Student t distributions. All distributions lie in the domains of attraction of GPDs with positive shape parameter $\xi_{DA} = 1, 0.25, 0.5, 0.25, 0.1$ and 0.05 respectively. In each case predictions are made from 10000 samples, and prediction performance is measured against a further 10000 samples. The predictor uses the first $k = 20$ upper order statistics, and the various lines correspond to sample sizes of $ N = 40 (.)$, $60 (+)$, $80 (*)$, $100 (\circ)$ and $200 (\Box)$. The $k=20$, $N = 200$ results are shown by the thicker lines.
  } \label{pred1}
\end{figure}

\begin{figure}[ph!] \centering
\includegraphics[height = 55mm,keepaspectratio]{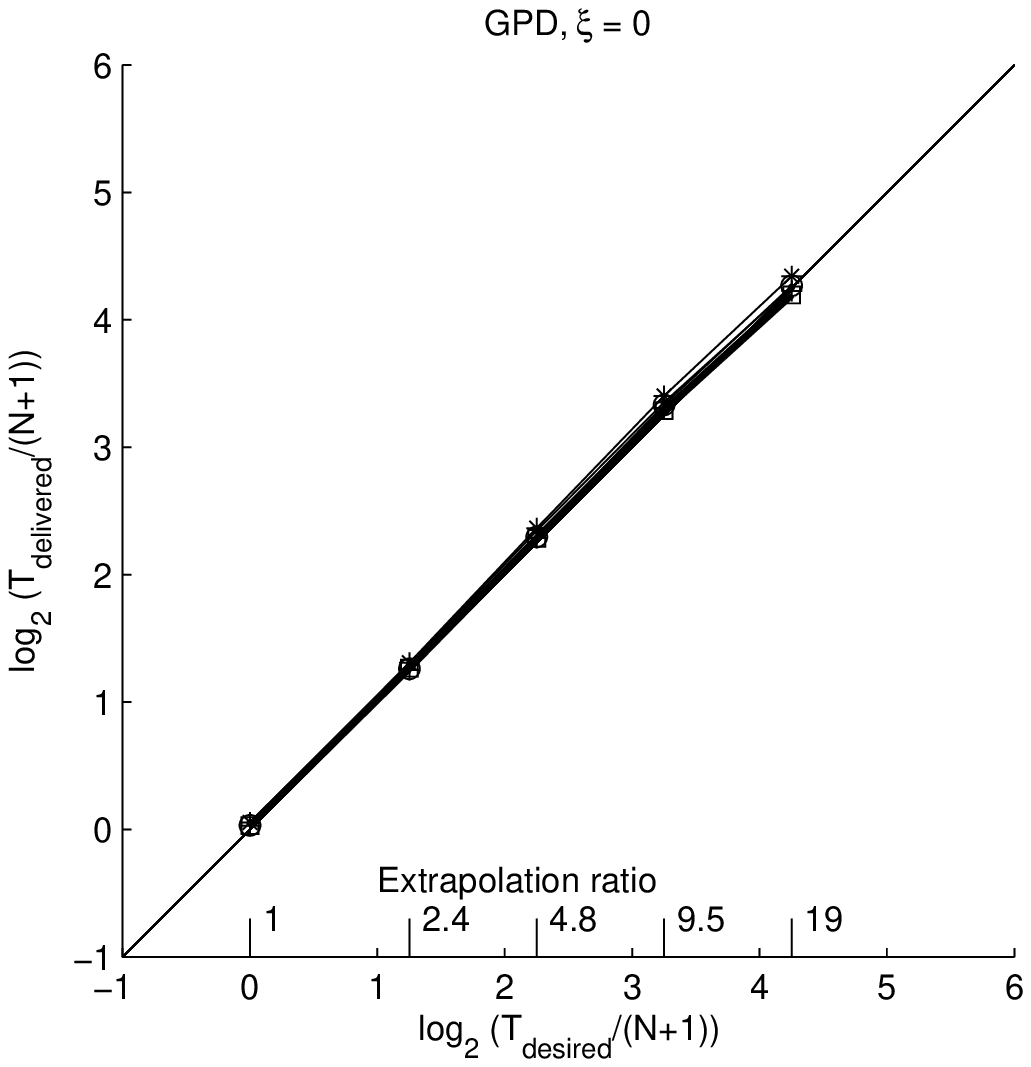}
\includegraphics[height = 55mm,keepaspectratio]{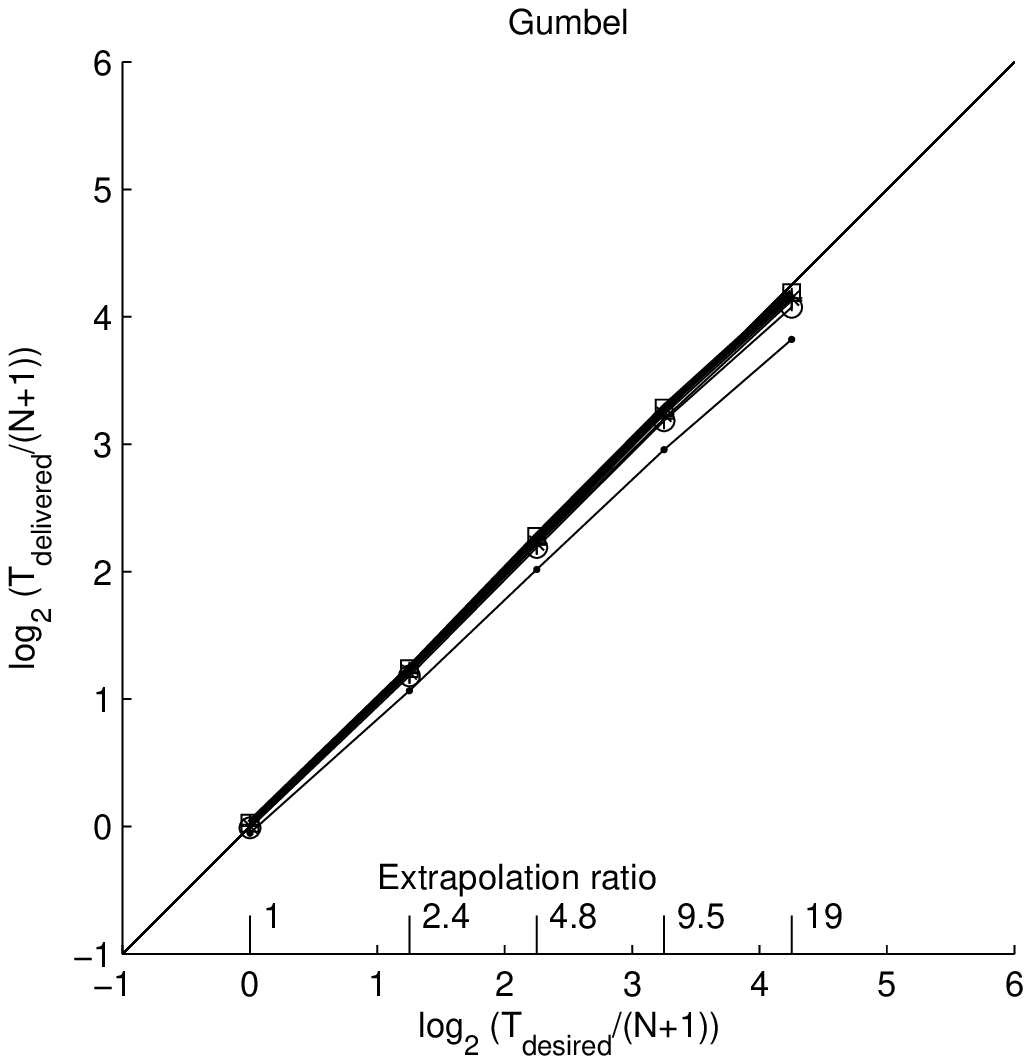} \\
\includegraphics[height = 55mm,keepaspectratio]{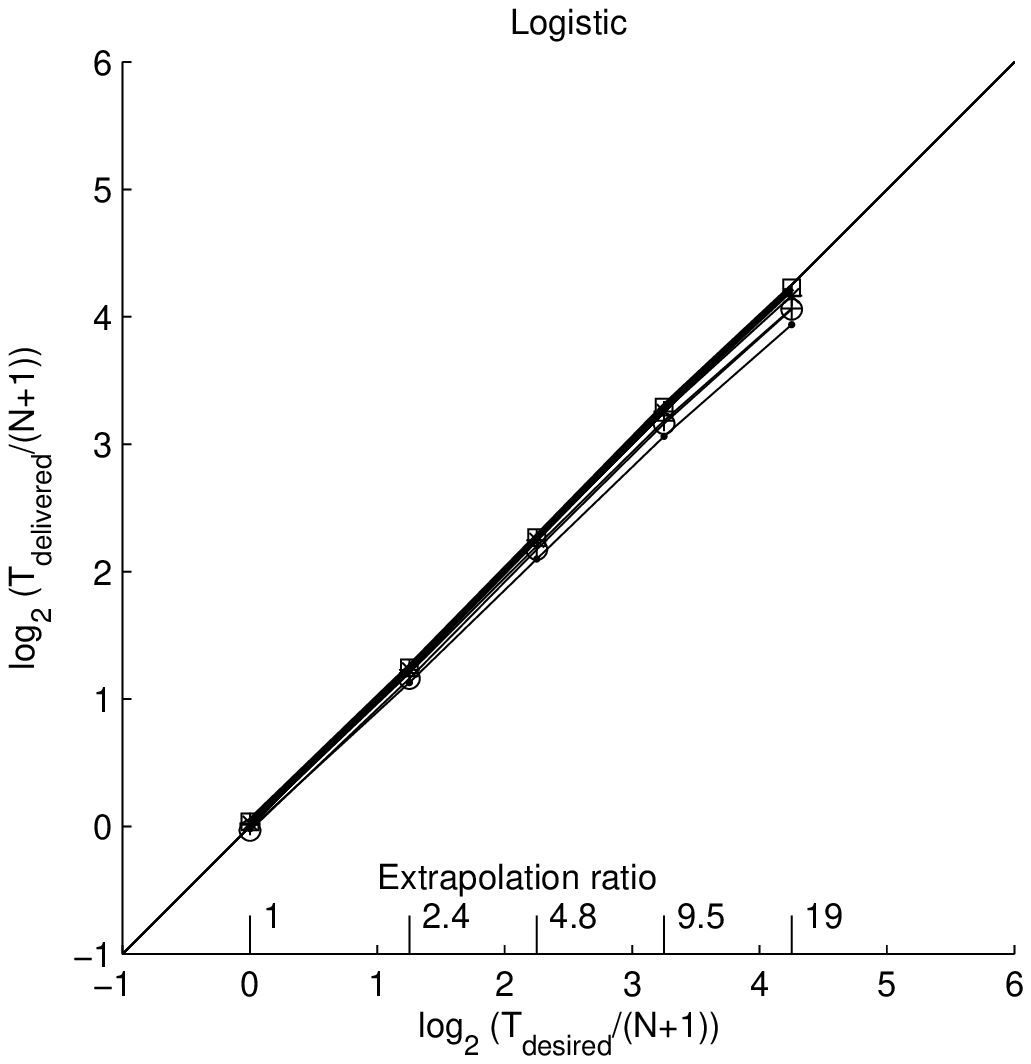}
\includegraphics[height = 55mm,keepaspectratio]{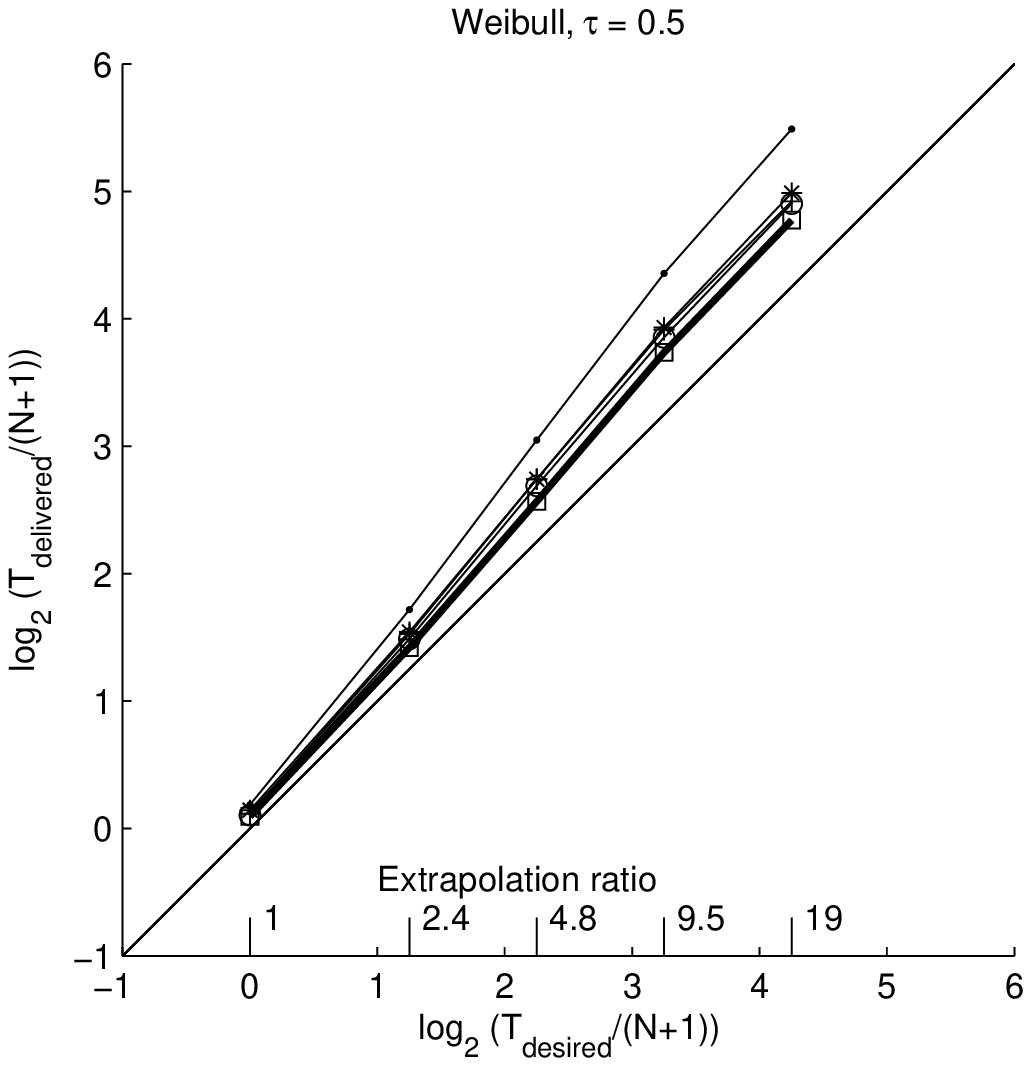} \\
\includegraphics[height = 55mm,keepaspectratio]{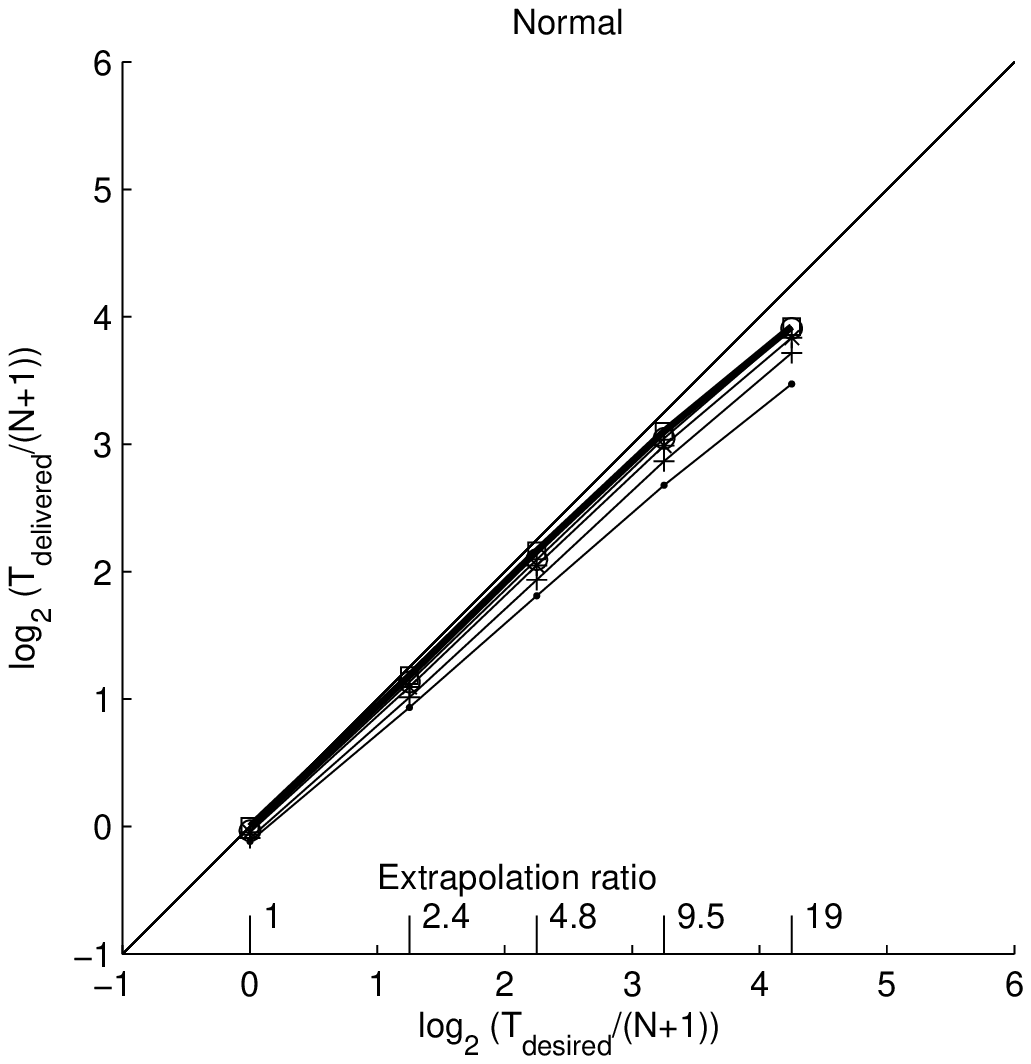}
\includegraphics[height = 55mm,keepaspectratio]{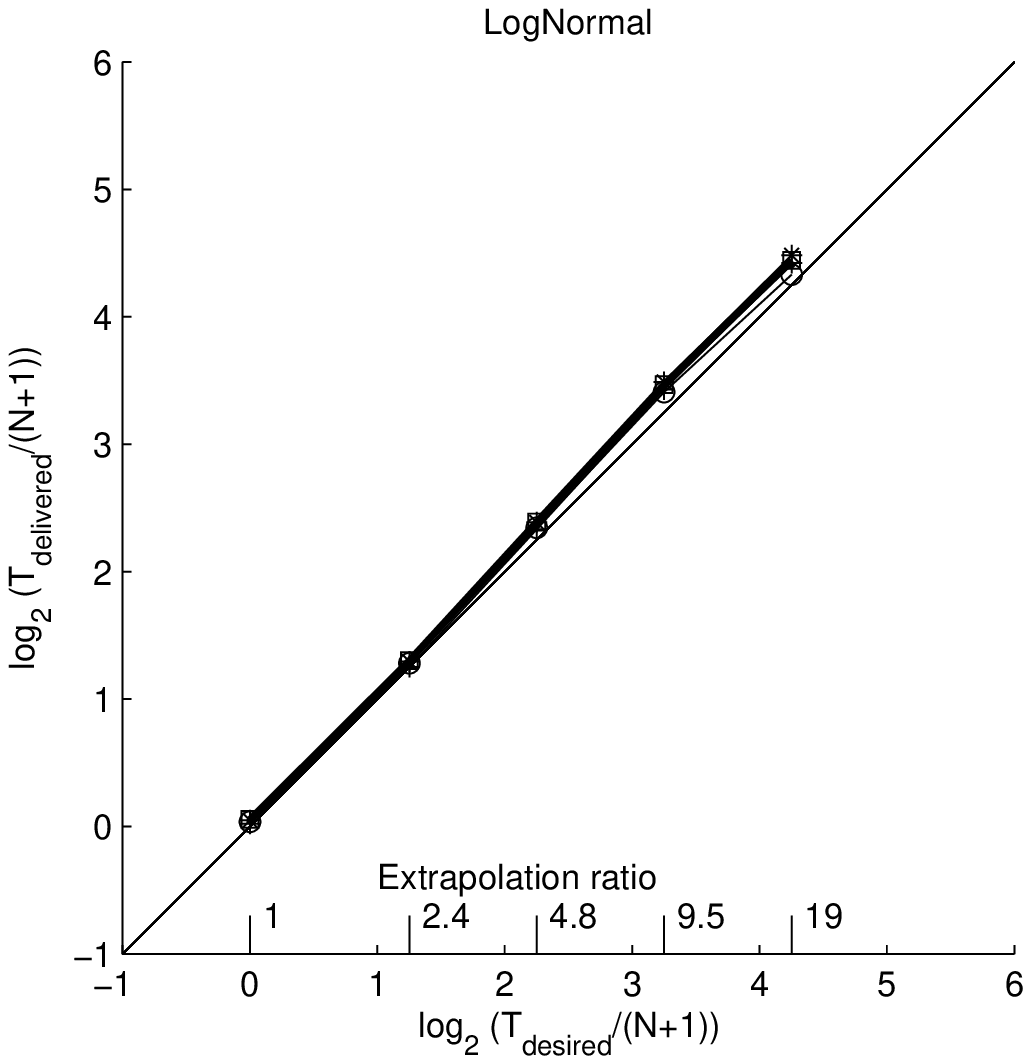}
  \caption{Probability performance of the predictor applied to samples drawn from the GPD ($\xi = 0$, Exponential), Gumbel Logistic, Weibull, Normal and Lognormal distributions. All distributions lie in the domain of attraction of the GPD with shape parameter $\xi_{DA} = 0$.
  }\label{pred2}
\end{figure}

\begin{figure}[ph!] \centering
\includegraphics[height = 55mm,keepaspectratio]{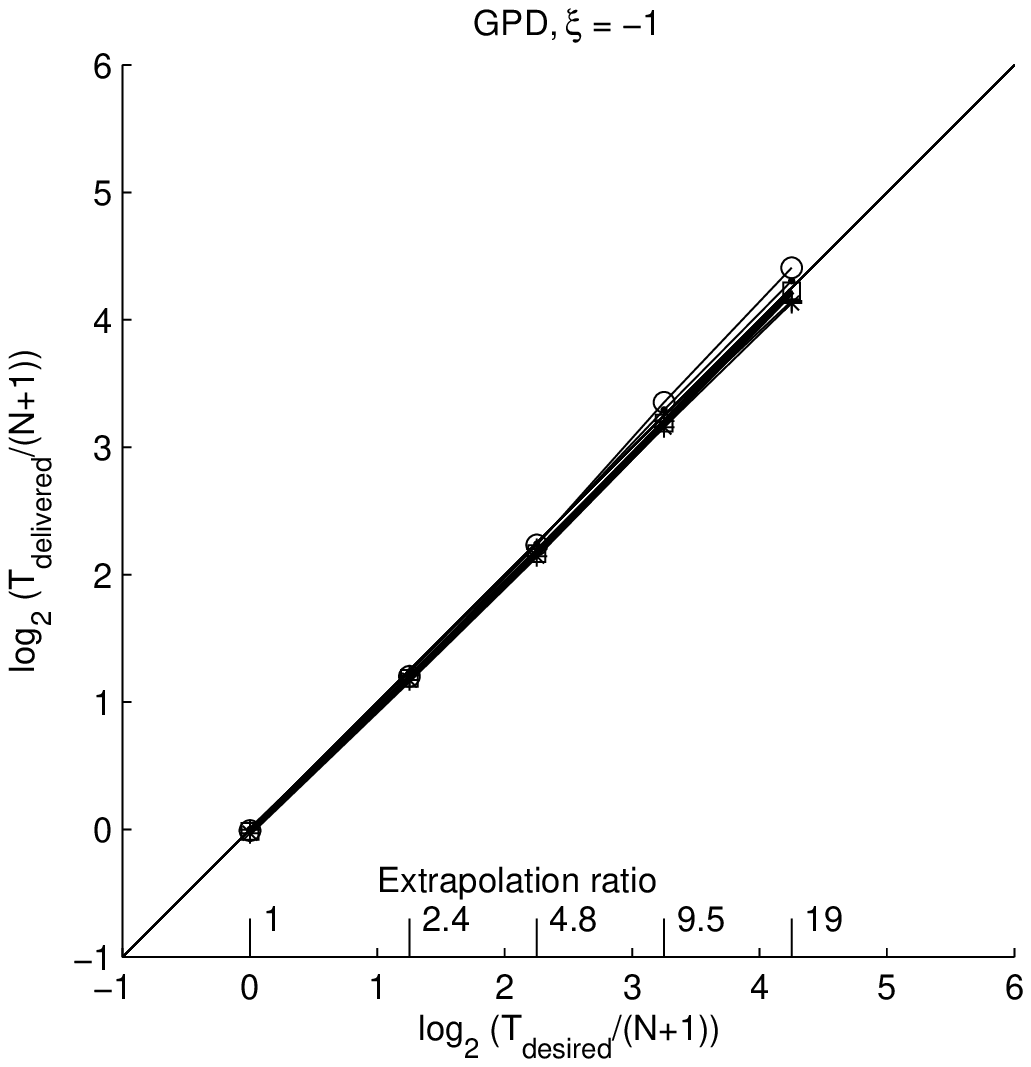}
\includegraphics[height = 55mm,keepaspectratio]{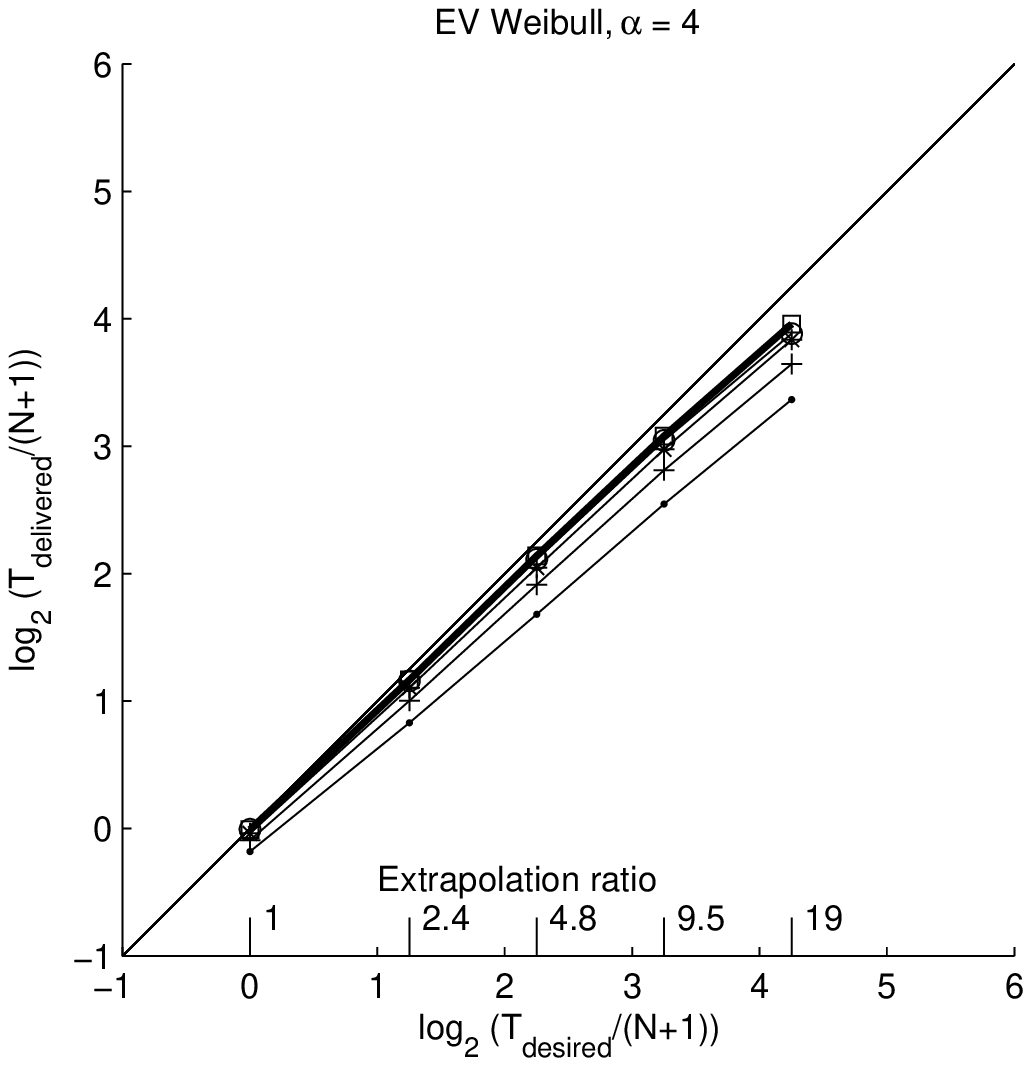} \\
\includegraphics[height = 55mm,keepaspectratio]{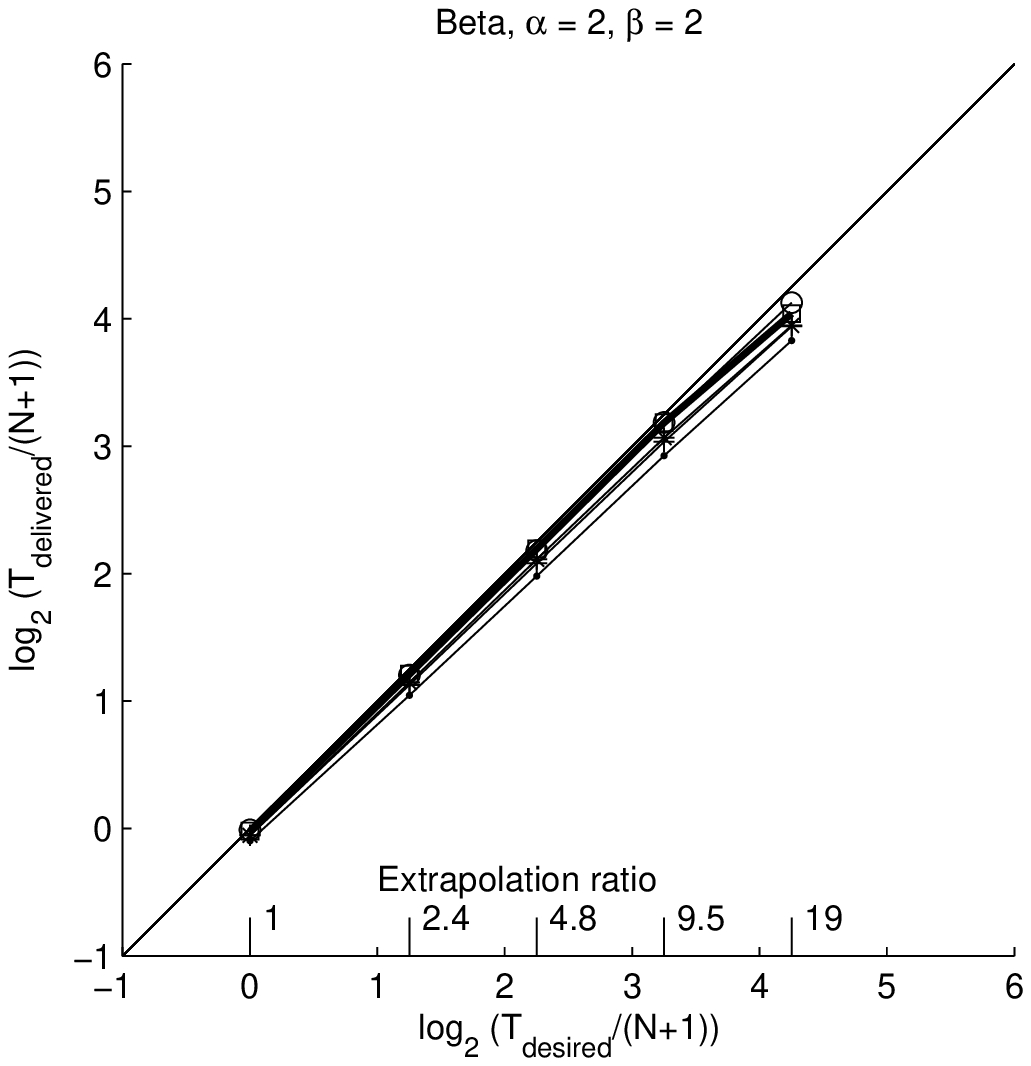}
\includegraphics[height = 55mm,keepaspectratio]{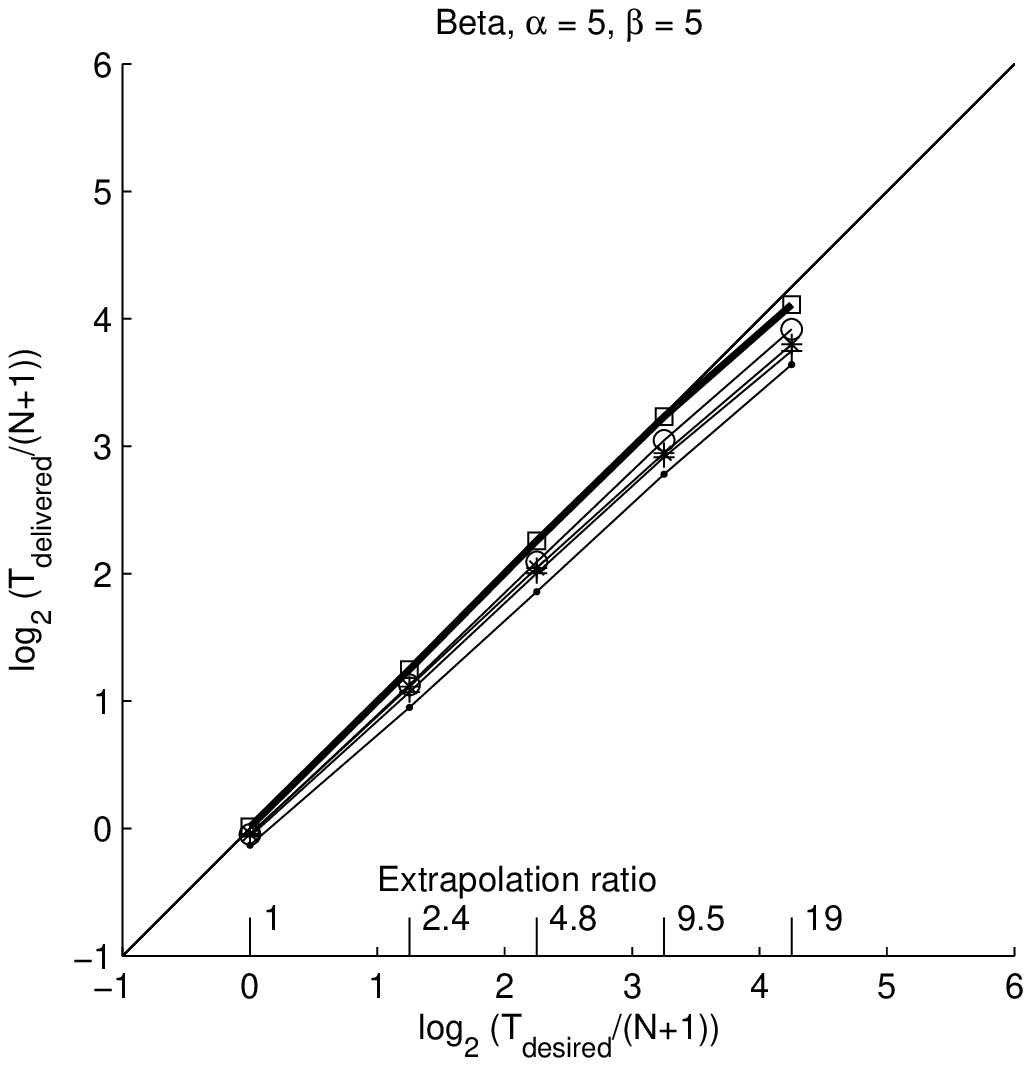} \\
\includegraphics[height = 55mm,keepaspectratio]{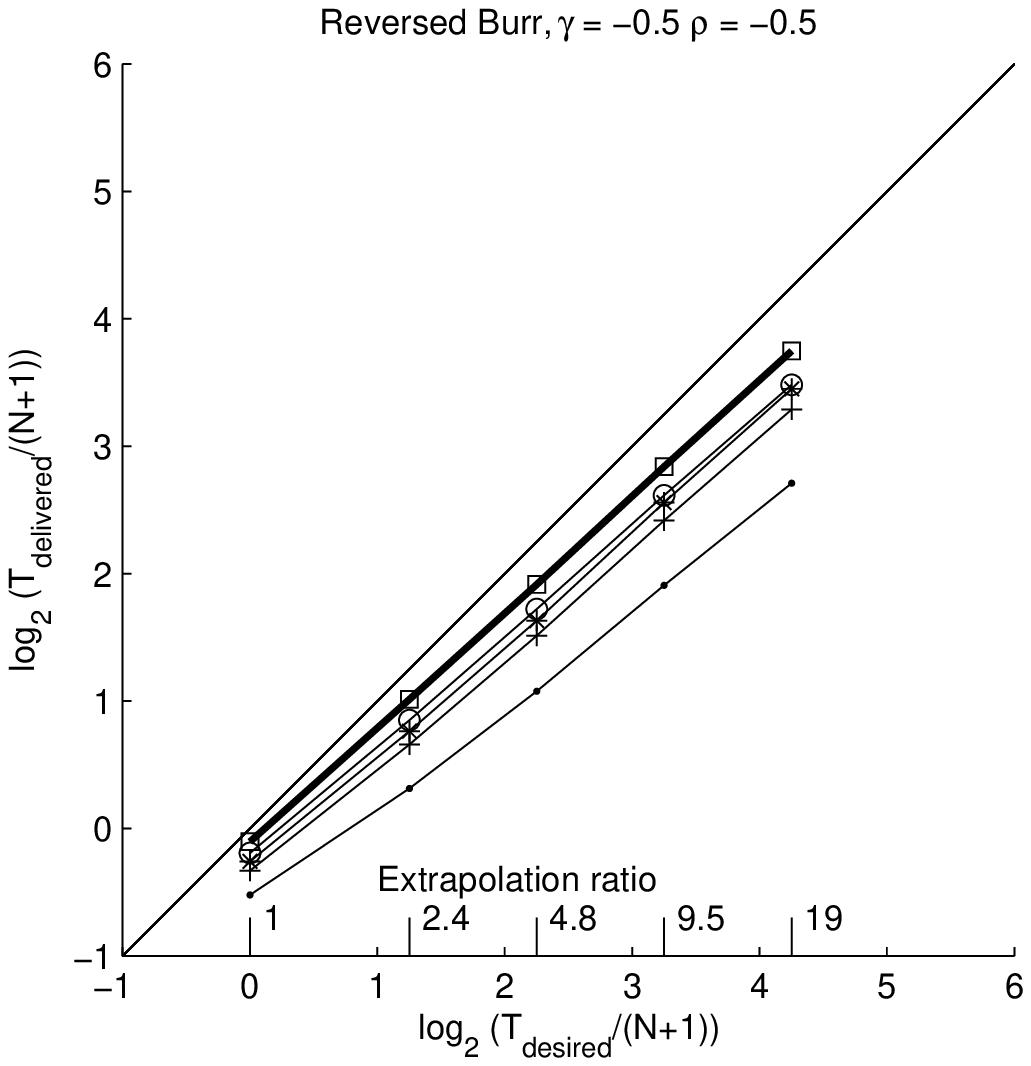}
\includegraphics[height = 55mm,keepaspectratio]{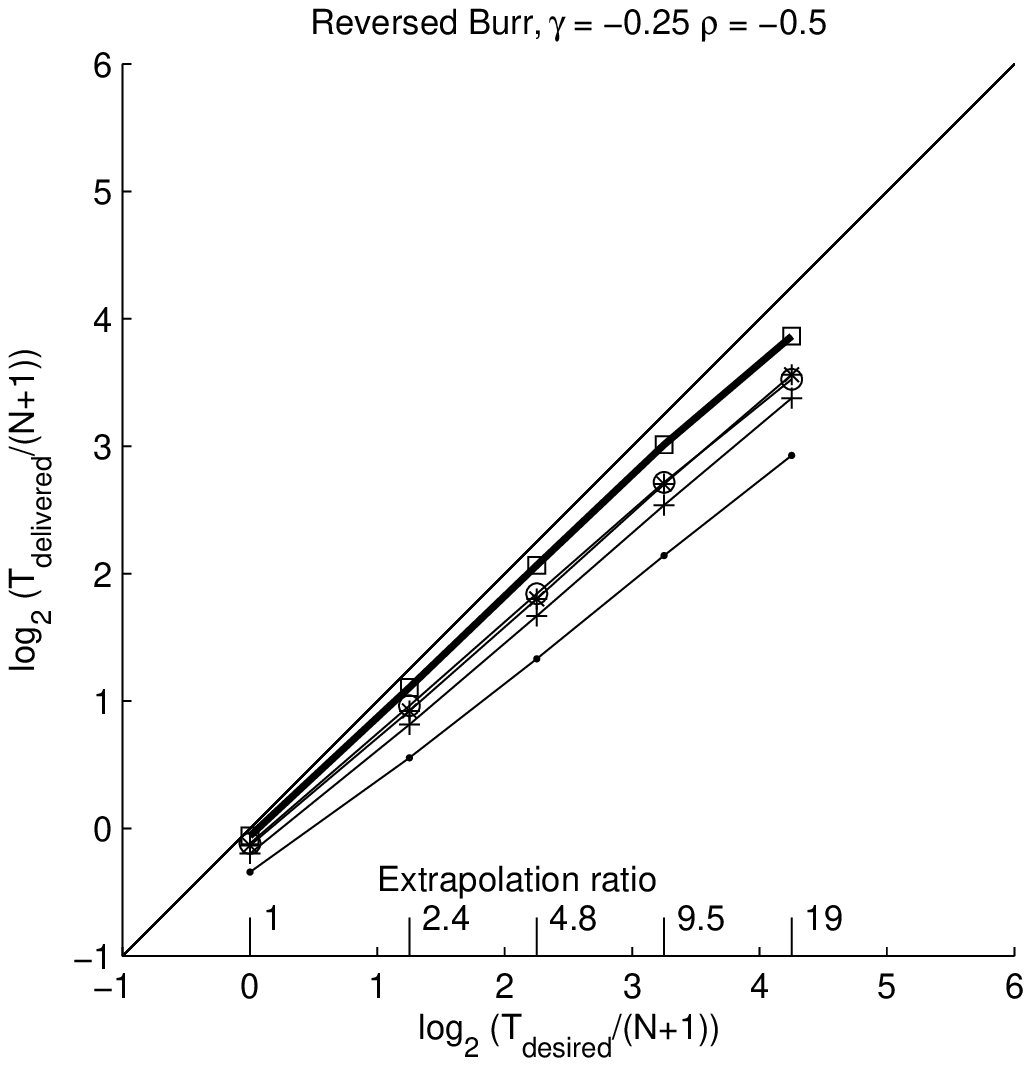}
  \caption{Probability performance of the predictor applied to samples drawn from GPD ($\xi = -1$), EV Weibull, Beta and Reversed Burr distributions. All distributions lie in the domains of attraction of GPDs with negative shape parameter $\xi_{DA} = -1, -0.25, -0.5, -0.2, -0.5$ and -0.25 respectively.
  }\label{pred3}
\end{figure}

For data drawn from a variety of  non-GPD distributions, Figures~\ref{pred1}-\ref{pred3} plots the exceedence level delivered compared with that desired. In each subfigure, five curves are shown, corresponding to increasing sample size through $N=40$, 60, 80, 100 to $N=200$, the latter results shown by the thicker lines. For $N = 40$ we do not expect to get good prediction performance since the predictor is constructed using half the data, which can scarcely be described as the  ``tail'', especially for double-sided distributions such as the normal. The key results are thus the thicker lines for $N=200$, where only the upper 10 percent of the data have been used for the extrapolation.

The non-GPD distributions considered correspond to those modelled in \cite{segers}. In all cases the delivered exceedance level is close to that desired, with the agreement improving for larger sample sizes. The largest extrapolation ratio shown is $400/21 \approx 19$ thus for the largest samples (of size $N=200$, shown by the thicker solid lines), the uppermost prediction is to a return level of $T = 3829$, which is over an order of magnitude beyond the return level ($T = N+1 = 201$) that would be naturally associated with the 200 data points. For such a large extrapolation, it is perhaps remarkable how closely the probability is preserved in all cases.

\pagebreak

\pagebreak

\section{Summary}
A method has been presented for making large extrapolations out-of-sample, yet in a way that gives good probability preservation.  It takes a data set of size $N \geq 20$, and uses only the highest 20 upper order statistics, whose distribution is then approximated by a member of the three parameter family of GPDs. Interest in the location and scale parameters is removed by rescaling with respect to the 10th and 20th upper order statistics. The shape parameter $\xi_{DA}$ of the domain of attraction GPD is then estimated as $\hat{\xi}$ using a least squares curve fit (which is accompanied by a useful graphical representation). Extrapolation at any desired return level $T_{des}$ is then accomplished using a GPD with an adjusted tail parameter $\xi_p = \hat{\xi}+d\xi$, where the increment $d\xi = d\xi(\hat{\xi},T_{des})$ is a function which has been calibrated to give good probability preservation for GPDs of any $\xi$. To a good approximation this probability preservation carries over to data drawn from a wide class of non-GPD distributions.

Only the basic outline of the method has been described and there is scope for further optimisation.

\bibliography{curvefitbib}

\begin{thebibliography}{7}
\expandafter\ifx\csname natexlab\endcsname\relax\def\natexlab#1{#1}\fi
\expandafter\ifx\csname url\endcsname\relax
  \def\url#1{\texttt{#1}}\fi
\expandafter\ifx\csname urlprefix\endcsname\relax\def\urlprefix{URL }\fi

\bibitem[{Castillo and Daoudi(2009)}]{castillo}
Castillo, J., Daoudi, J., 2009. Estimation of the generalized {P}areto
  distribution. Statistics and Probability Letters 79, 684--688.

\bibitem[{Embrechts et~al.(1999)Embrechts, Kl\"{u}ppelberg, and
  Mikosch}]{embrechts}
Embrechts, P., Kl\"{u}ppelberg, C., Mikosch, T., 1999. Modelling Extreme Events
  for Insurance and Finance. Springer, Berlin.

\bibitem[{Mathworks(2014)}]{gpfit}
Mathworks, 2014. Matlab function {\it gpfit}.

\bibitem[{McRobie(2013{\natexlab{a}})}]{McRobieGEV}
McRobie, F.~A., 2013{\natexlab{a}}. Elemental estimators for the {G}eneralized
  {E}xtreme {V}alue tail. ar{X}iv:1304.4362.

\bibitem[{McRobie(2013{\natexlab{b}})}]{McRobieGPD}
McRobie, F.~A., 2013{\natexlab{b}}. Elemental unbiased estimators for the
  {G}eneralized {P}areto tail. ar{X}iv:1304.3918.

\bibitem[{McRobie(2013{\natexlab{c}})}]{McRobiePRED}
McRobie, F.~A., 2013{\natexlab{c}}. Probability-matching predictors for extreme
  extremes. ar{X}iv:1307.7682.

\bibitem[{Segers(2005)}]{segers}
Segers, J., 2005. Generalized {P}ickands estimators for the extreme value
  index. J. Stat. Planning and Inference 128~(2), 381--396.

\end{thebibliography}
\section{Appendix: Insight into tail estimation errors for non-GPD distributions}

We demonstrate briefly how the basic curve-fit construction can give insight into cases where estimators for the tail index of  non-GPD distributions lead to persistently erroneous results.

\begin{figure}[h!] \centering
\includegraphics[width=48mm,keepaspectratio]{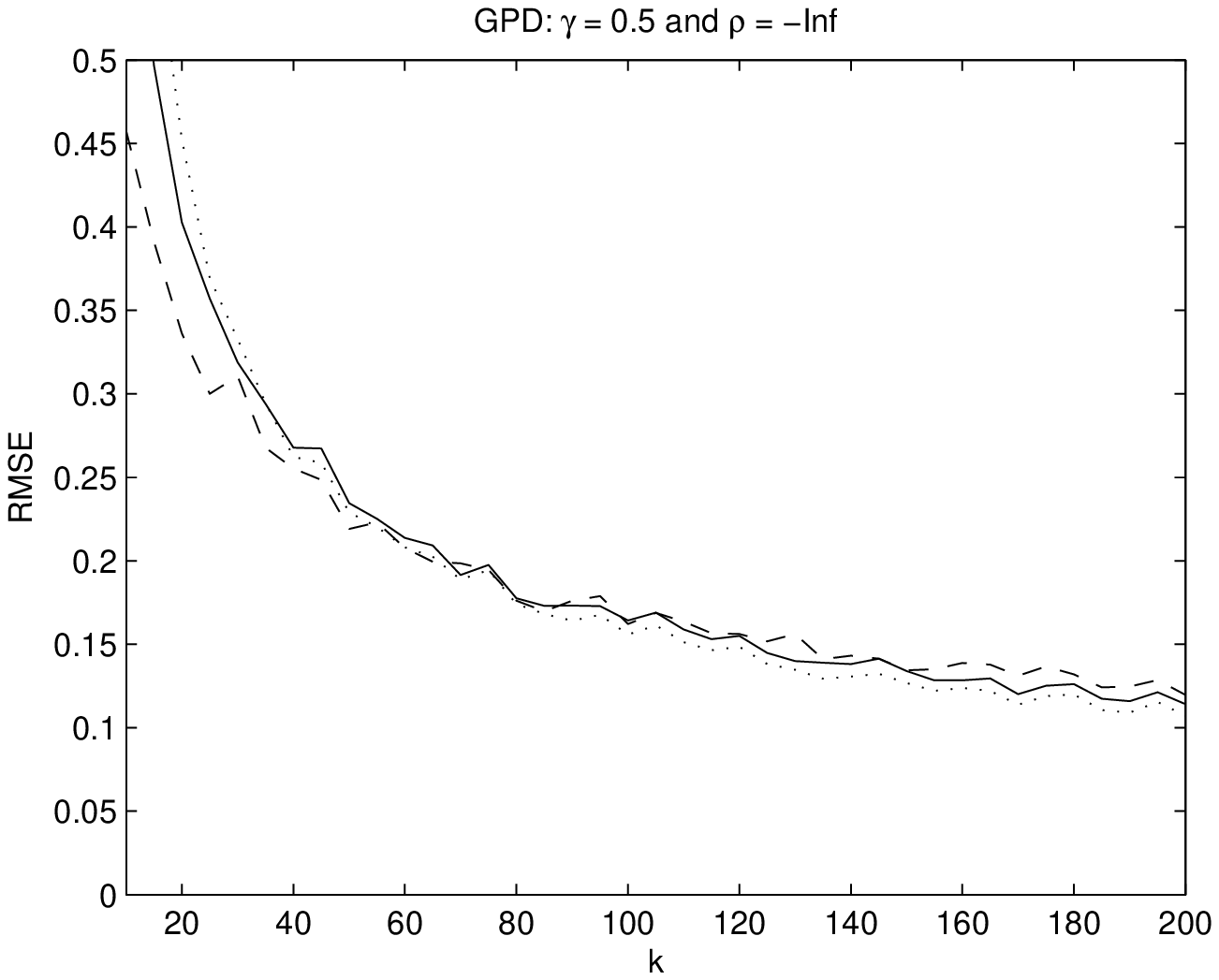}
  \includegraphics[width=48mm,keepaspectratio]{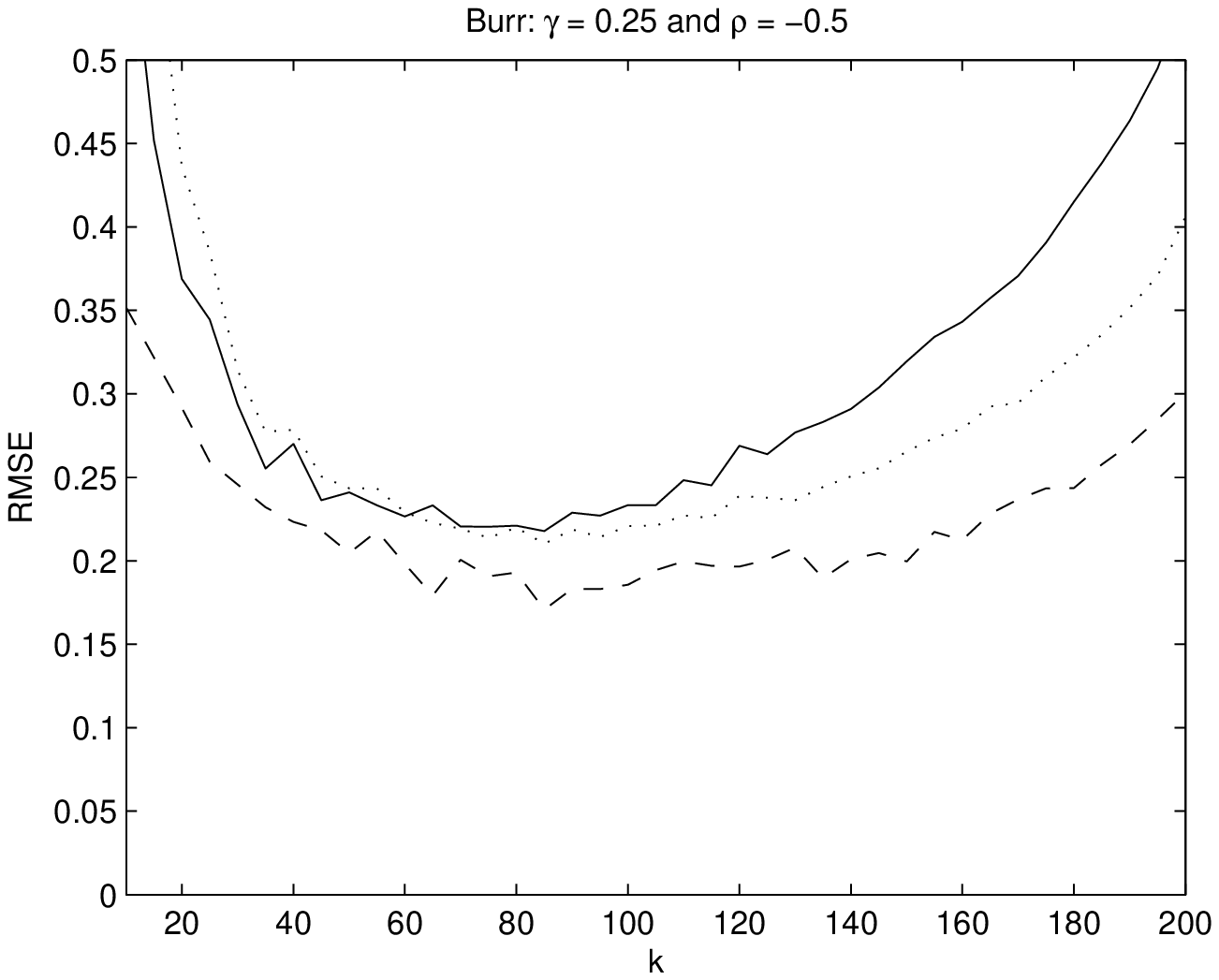}
  \includegraphics[width=48mm,keepaspectratio]{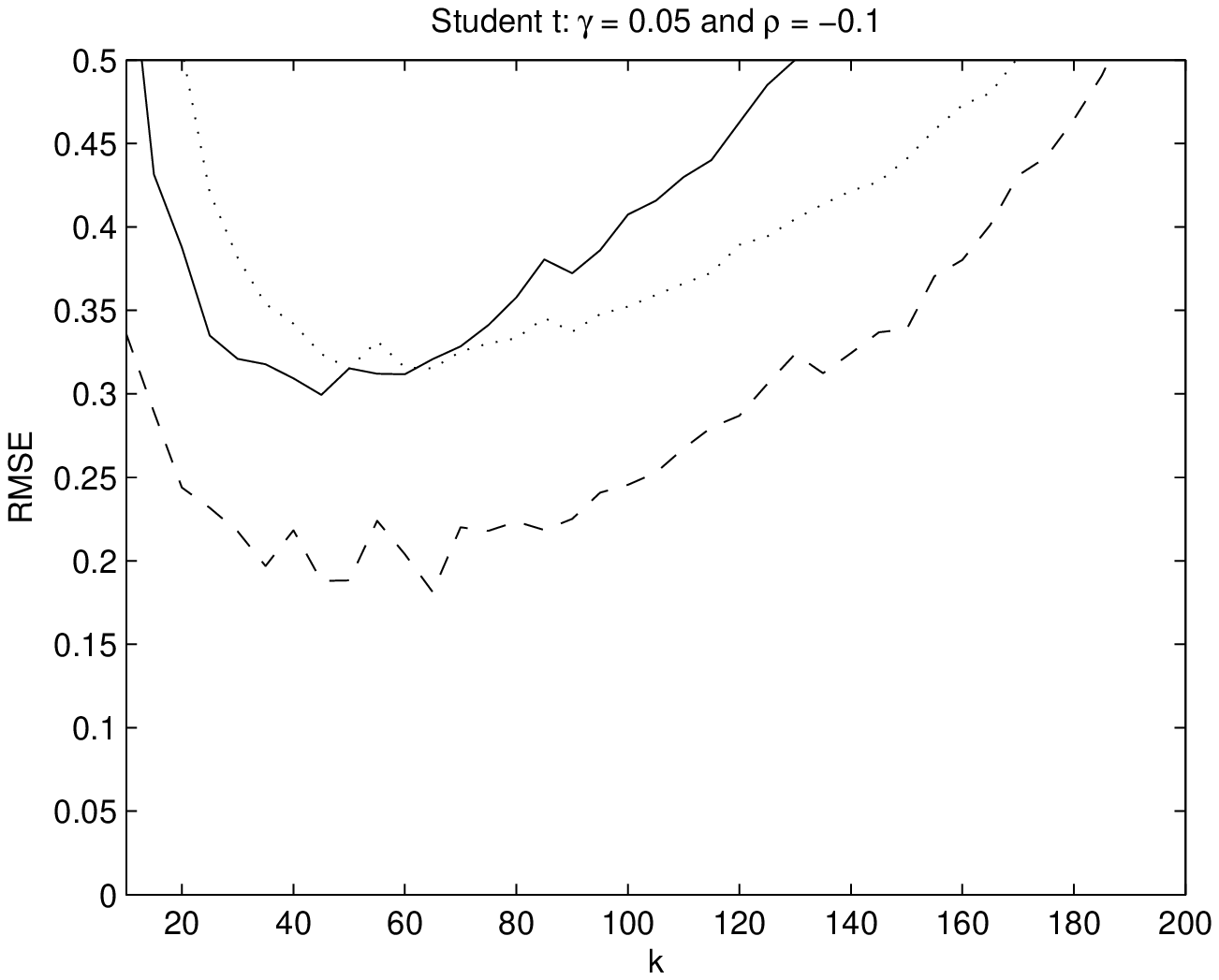}\\
  \includegraphics[width=48mm,keepaspectratio]{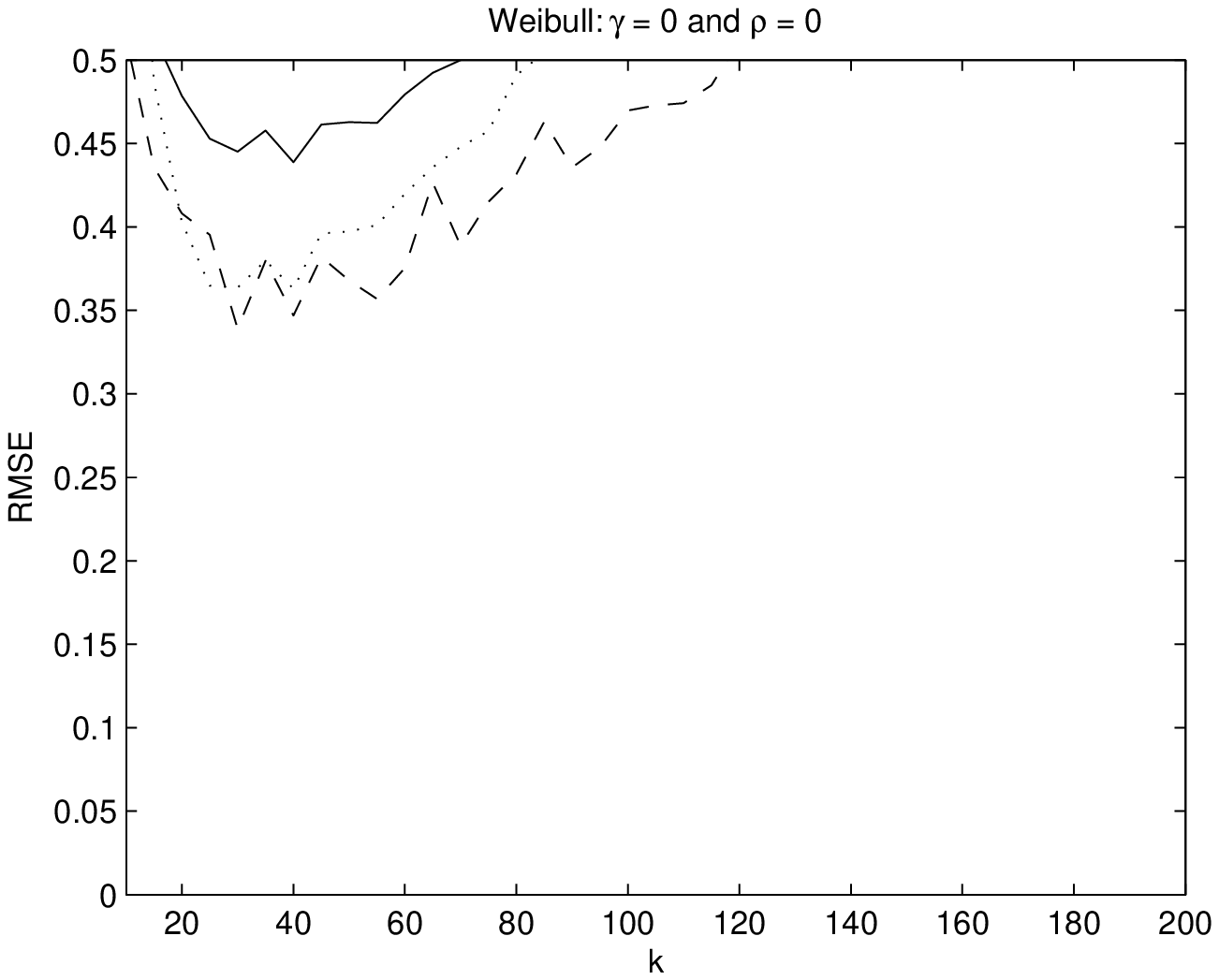}
  \includegraphics[width=48mm,keepaspectratio]{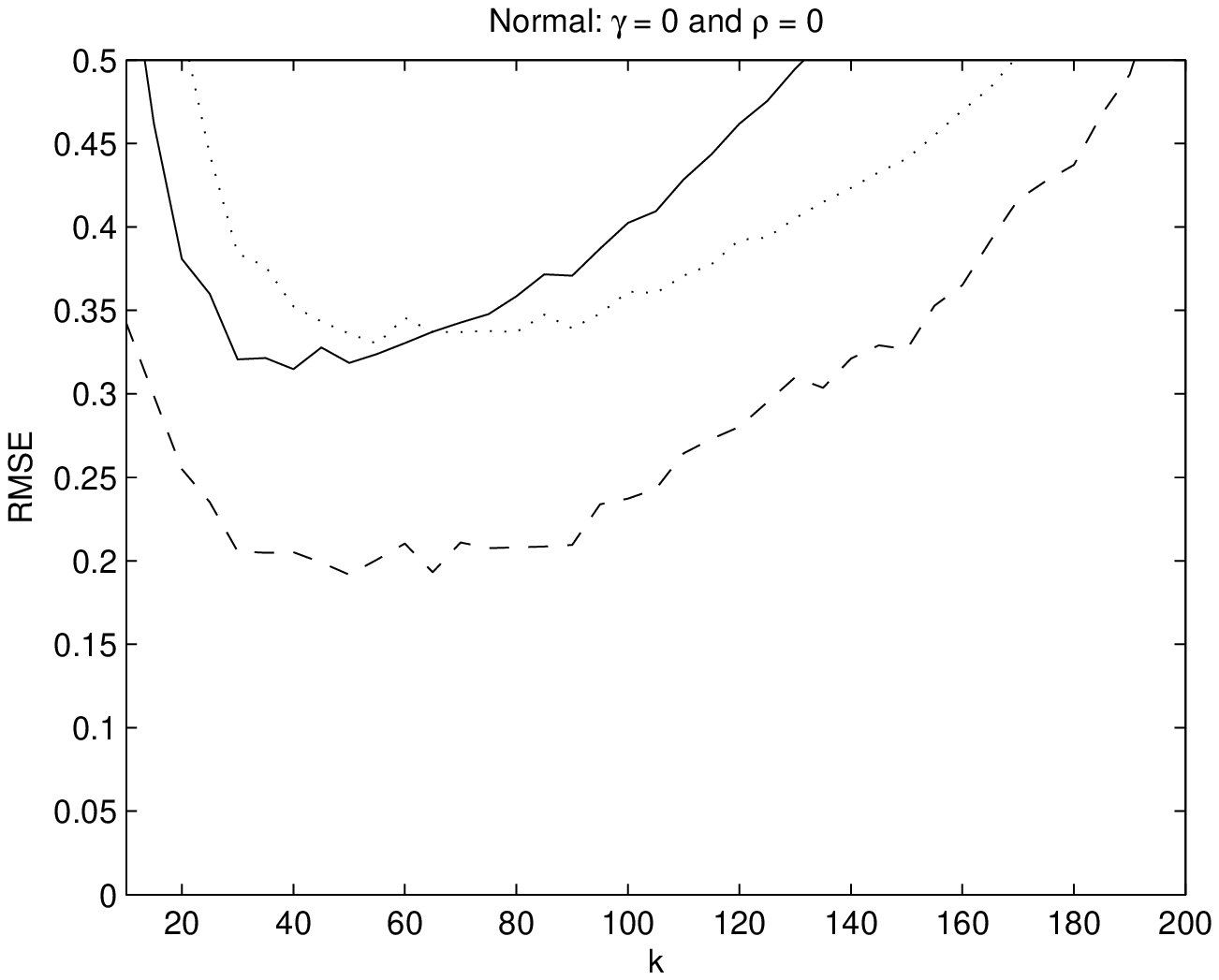}
  \includegraphics[width=48mm,keepaspectratio]{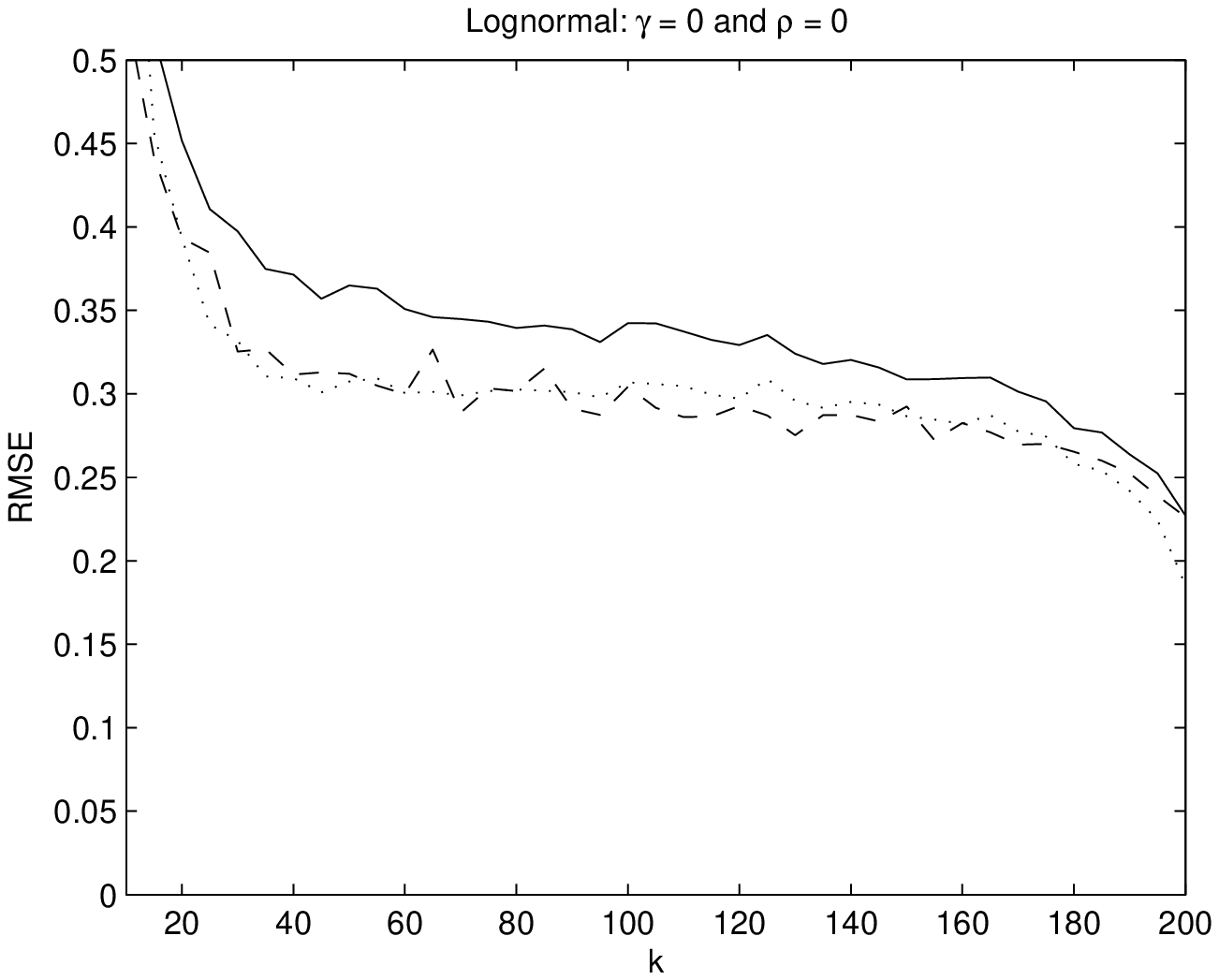}\\
  \includegraphics[width=48mm,keepaspectratio]{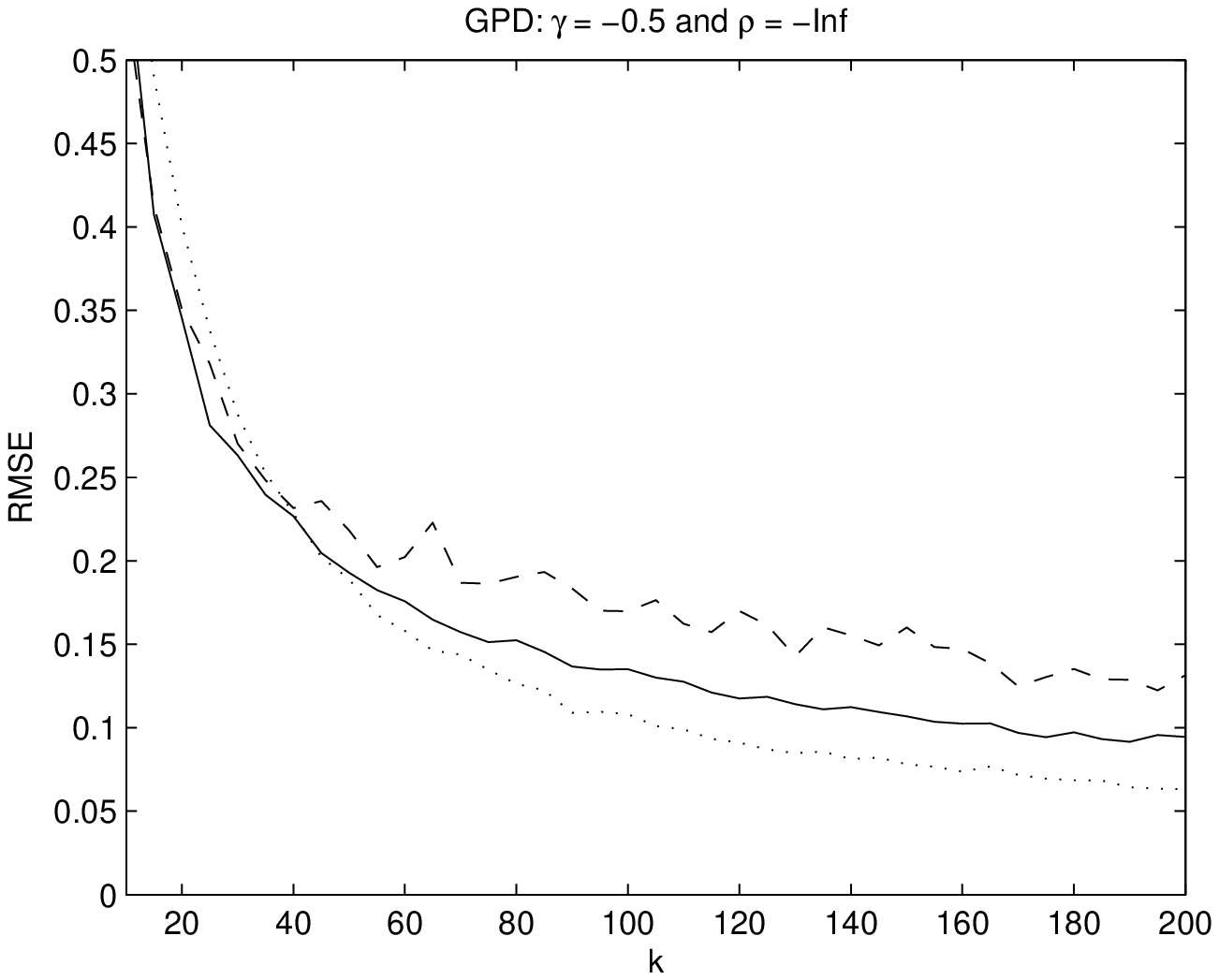}
  \includegraphics[width=48mm,keepaspectratio]{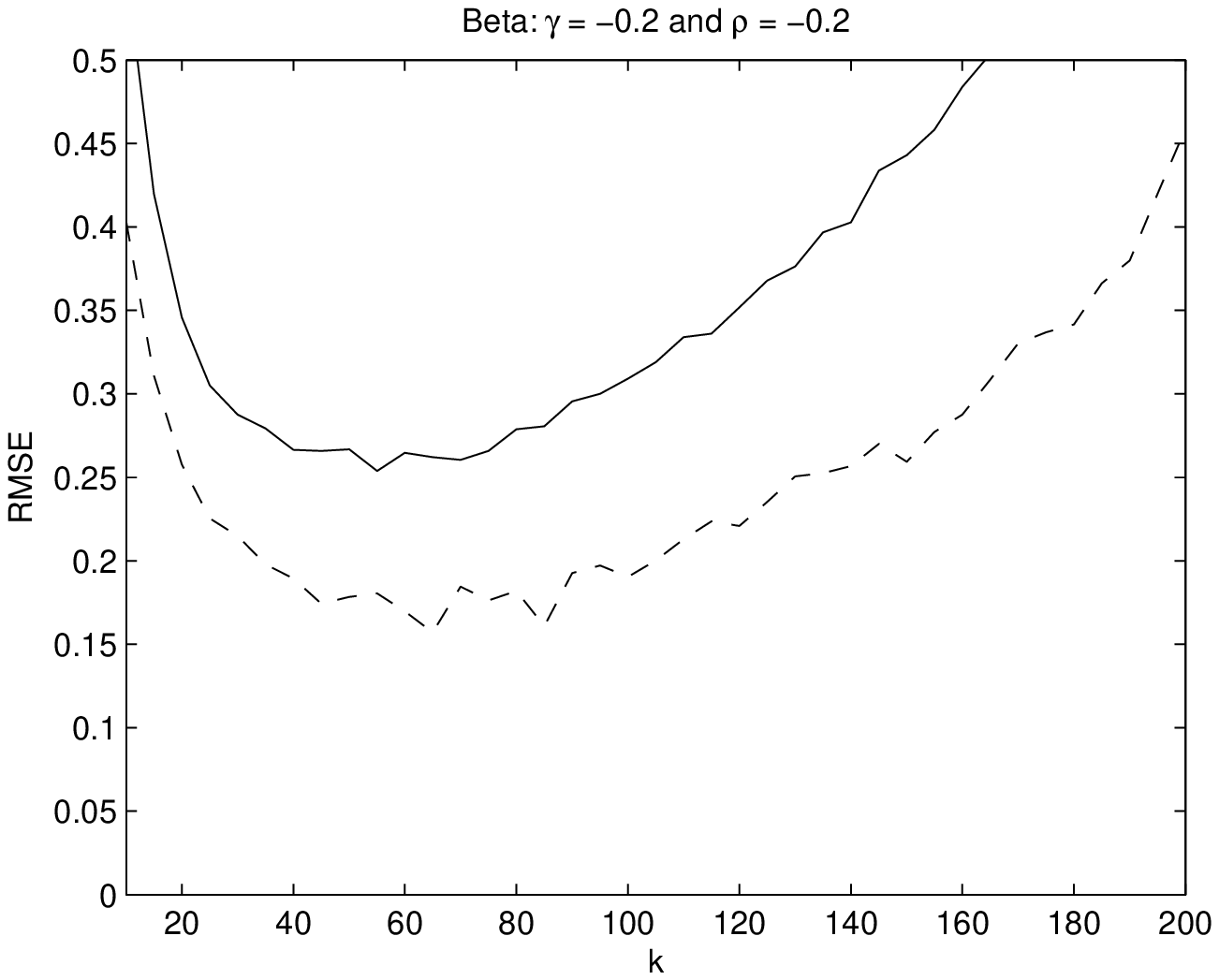}
  \includegraphics[width=48mm,keepaspectratio]{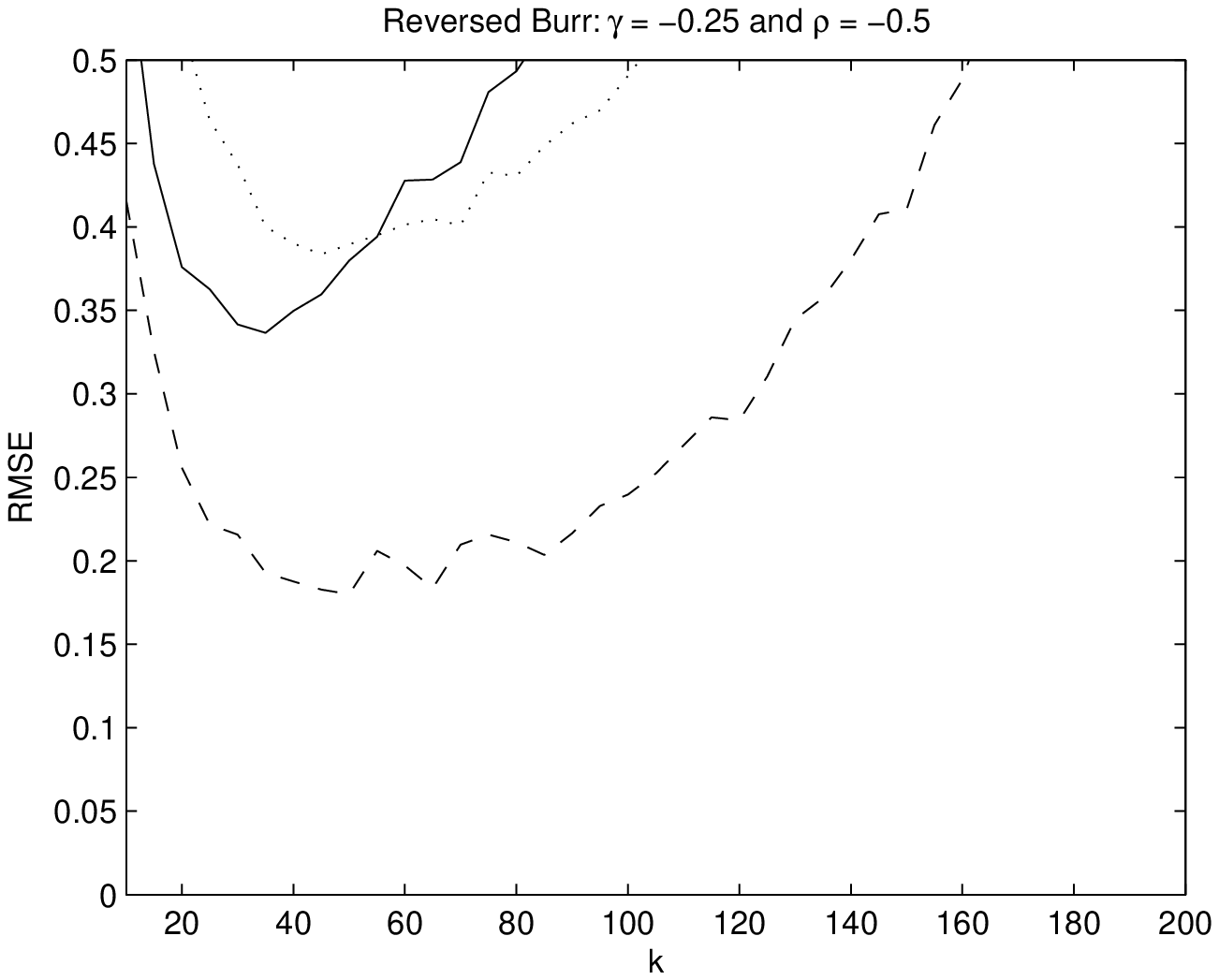}
  \caption{The root mean square error
  for the combined elemental (solid), the unconstrained Segers (dashed)
  and the maximum likelihood (dotted) estimators. The top, middle and bottom
  rows are for distributions in the domain of attraction of a GPD
  with positive, zero and negative tail parameters respectively. The
  nine cases (GPD, Burr, Student; Weibull, normal, lognormal;  GPD, beta, reversed Burr)
  are identical to ones used in \cite{segers}, where more
  detailed specification may be found.
  }\label{genfig1}
\end{figure}

Figure~\ref{genfig1} compares the root mean square errors (RMSE) of three tail index estimators applied to the $k$ upper order statistics of samples of size $N= 200$ from a wide range of distributions. The three estimators are the unconstrained Segers estimator (\cite{segers}), the linearly rising combination of elemental estimators (\cite{McRobieGPD}) and the Maximum Likelihood estimate of the Matlab {\it gpfit} function (\cite{gpfit}). As $k$ increases, the RMSE in almost the all non-GPD cases shows initial improvements before diverging again as $k$ approaches $N$.
In almost all cases, the Segers estimator has RMSE lower than the other two (the $\xi = -0.5$ GPD (lower left) being one exception). However, a focus on RMSE alone does not explain what is happening here. We focus on the case centre left, the $\tau = 0.5$ Weibull distribution, this being an example where all estimators are performing rather poorly.

\begin{figure}[h!] \centering
\includegraphics[height = 70mm,keepaspectratio]{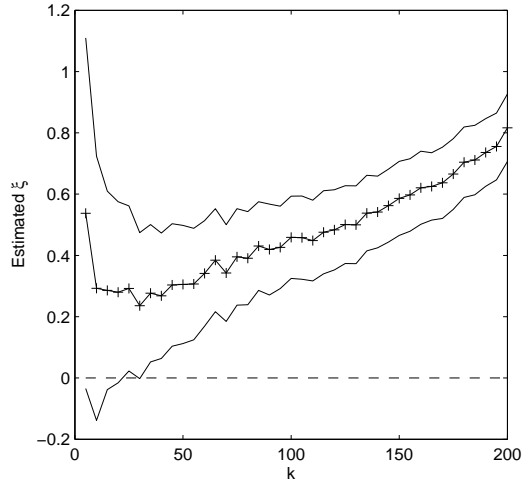}
  \caption{The unconstrained Segers estimator applied to 1000 samples of size $N = 200$ drawn from a Weibull distribution with parameter $\tau = 0.5$.  The estimates of the tail index $\xi_{DA}$ are plotted as ordinate. The estimator is constructed using only the first $k$ upper order statistics, with $k$ plotted as abscissa.
  }\label{segersweib}
\end{figure}

Figure~\ref{segersweib} shows the unconstrained Segers tail index estimates for the $\tau= 0.5$ Weibull, plotting mean $\pm$ standard deviation. The theoretical answer for a Weibull is $\xi_{DA} = 0$. The estimator has been applied to the first $k$ upper order statistics of samples of size $N=200$. (The corresponding plot for the other two estimators are decidedly similar, with estimates lying just above those of the Segers estimator.) It can be seen that, as $k$ increases,  the reductions in RMSE that are associated with the reduced variance of working with more data points are soon outweighed by the increased bias as $k$ approaches the full sample size $N$.

However, Figure~\ref{segersweib} suggests that the question of which estimator has the marginally better RMSE is perhaps of lesser interest than the question as to why all the estimates are so far from the correct value. Although the true answer is $\xi_{DA} = 0$, all the estimators seems fairly convinced that $\xi_{DA}$ is actually positive.

That an extreme value estimator performs badly for the Weibull distribution is perhaps surprising, given the comparative ubiquity of Weibulls in the extreme value literature. Weibulls are commonly applied in reliability analysis, such as to the failure rates of components. Also, there is a form of the Weibull inherent in the family of Generalised Extreme Value distributions, although this is the reversed- or EV-Weibull (for which tail index estimators work rather well).

 Some insight into why the estimators perform so badly on the (unreversed) Weibull can be obtained by plotting Weibull information onto the graphical construction of the main body of this paper.

The distribution function of the three-parameter Weibull may be written
\begin{equation}
F(x) = 1 - \exp \left( - \left( \frac{x-\mu}{\sigma} \right)^{\tau} \right)
\end{equation}

\begin{figure}[h!] \centering
\includegraphics[height = 75mm,keepaspectratio]{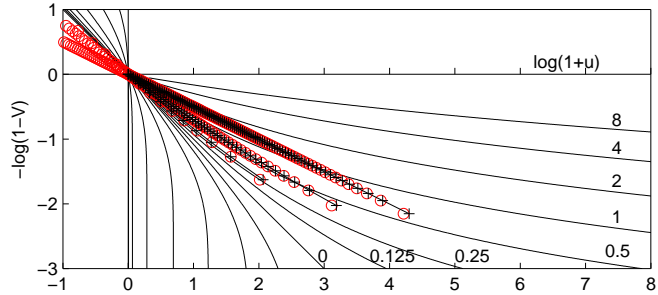}
  \caption{ Weibull samples of size $N=200$ plotted on the intuitive curve-fit tail plot (using the $G_i = (i-0.5)/N$ plotting positions throughout).  From left to right, the three lines correspond to $k = 20$, 100, 200. Analytical estimates are denoted (+, black) and averages of 10000 samples are plotted (o, red).
  }\label{weibappenplot}
\end{figure}

Inversion of the tail distribution function, as was done for the GPD, leads to the analytical approximation
\begin{equation}
\tilde{u}_i = \frac{g_i^{1/\tau} - 1}{1 - \alpha^{1/\tau}} \text{   with }  g_i = \frac{\log G_i}{\log G_j} \text{  and  } \alpha = \frac{\log G_k}{\log G_j}
\end{equation}
In Figure~\ref{weibappenplot}, these analytical approximations are plotted onto the $k- 0.5$  construction for $N= 200$ for the three cases $k = 20,$ 100 and 200.

The insight afforded by the construction is now apparent. The plot makes visible the way that the tail is ``pulled in'' towards the true solution $\xi_{DA} = 0$  as $k$ is decreased. Even by eye, the picture suggests that a GPD-based tail estimator would be expected to fall from a value just below 1.0 for $k= N = 200$ to a value near 0.25  for $k= 20$. Applying least squares to the analytical approximations gives
a slightly more precise view of what should be anticipated when a GPD predictor is applied to Weibull data: the tail estimates should be expected to be around  $\hat{\xi} = 0.21$, 0.51 and 0.93 for $k=20$, 100 and 200 respectively. The performance of the unconstrained Segers estimator (Figure~\ref{segersweib}) is decidedly similar.

Averages of data sampled from Weibulls have also been plotted (in red) onto Figure~\ref{weibappenplot}, and these lie close to the analytical approximations. The analytical approximations, though, are the important information here, having been made {\it a priori}, without sampling. If the analytical form of a non-GPD distribution is known, then the basic construction allows the erroneous estimates that a GPD-based estimator may make to be largely anticipated {\it a priori} in a highly visual and intuitive manner.

\subsection{Corollary - a tail estimator for non-GPD data}

Given that the basic construction shows how tail estimates may change under increasing $k$ for non-GPD data, it suggests that such changes may be monitored to assess just how non-GPD the data is, leading to the possibility of improved estimates.

For example, at any $k$ one could make a sequence of curve fit estimates $\hat{\xi}'$ at a sequence of $k'$ which increase to $k$ and look at how the tail is pulled in as $k'$ increases. Some careful back-extrapolation of the $\hat{\xi}'$ back to the hypothetical value $k' = 0$ may thus lead to reduced bias in the level $k$ estimate.

Figure~\ref{weibappenplot} gives the visual insight. There, the sequence of progressively smaller $k$ leads to a sequence of curves that approach the true value $\xi_{DA} = 0$. When using $k=200$, is the analyst supposed to accept an estimate $\hat{\xi} \approx 1$, even though those 200 data points contain within them the $k' = 20$ data which is suggesting $\hat{\xi} \approx 0.25$?
 Thus when using $k = 200$, rather than merely accepting the direct $k = 200$ estimate of $\hat{\xi} \approx 1$ or even the $k' = 20$ estimate of $\hat{\xi} \approx 0.25$, a value lower than 0.25 appears to be reasonable. This could be estimated by back-extrapolation of the various estimates $\hat{\xi}'$ back to the hypothetical $k' = 0$ intercept.  Of course, Figure~\ref{weibappenplot} shows averages, whereas the procedure would need to be applied to individual samples and the back-extrapolation results would be complicated by the associated scatter.

\begin{figure}[h!] \centering
\includegraphics[height = 50mm,keepaspectratio]{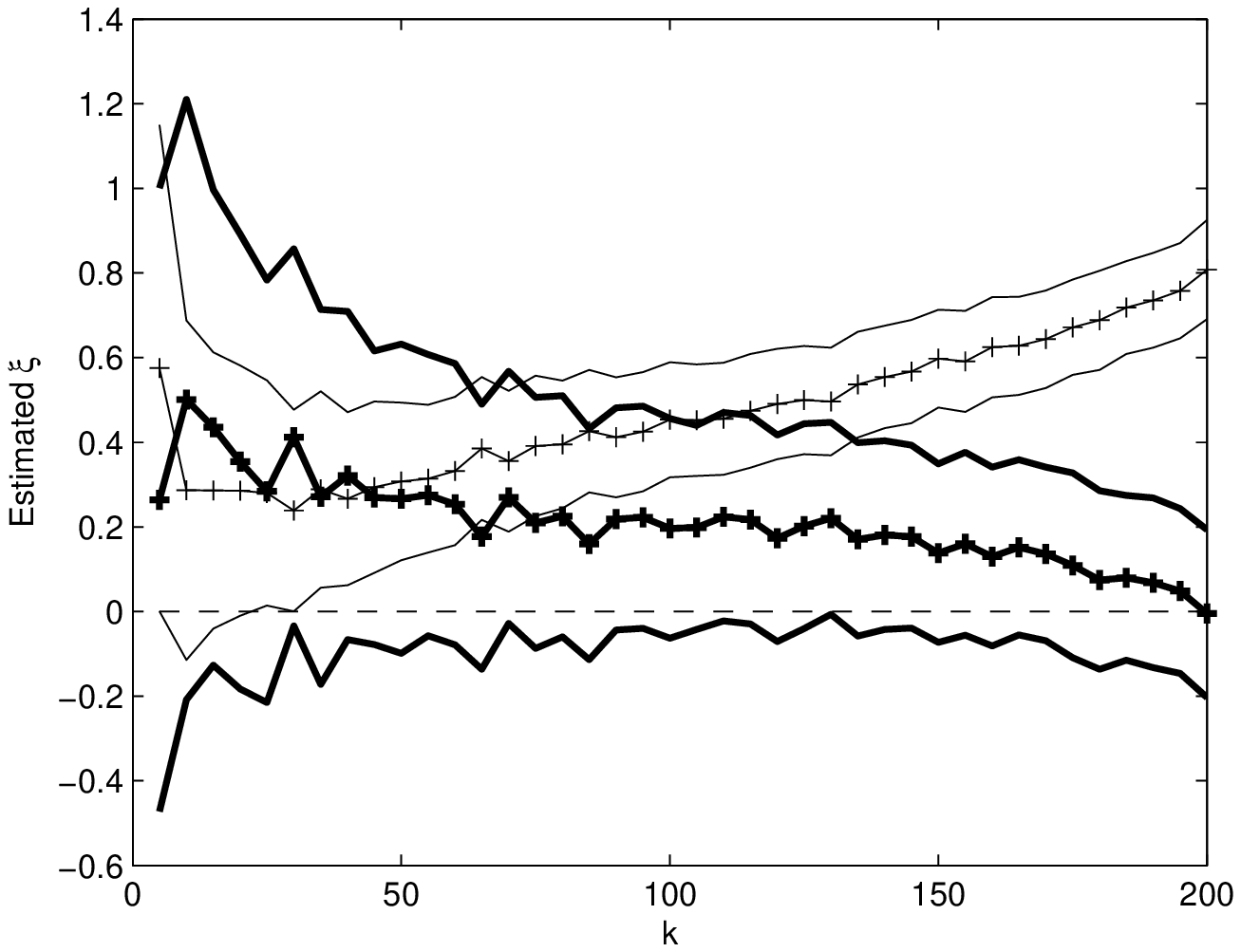}
\includegraphics[height = 50mm,keepaspectratio]{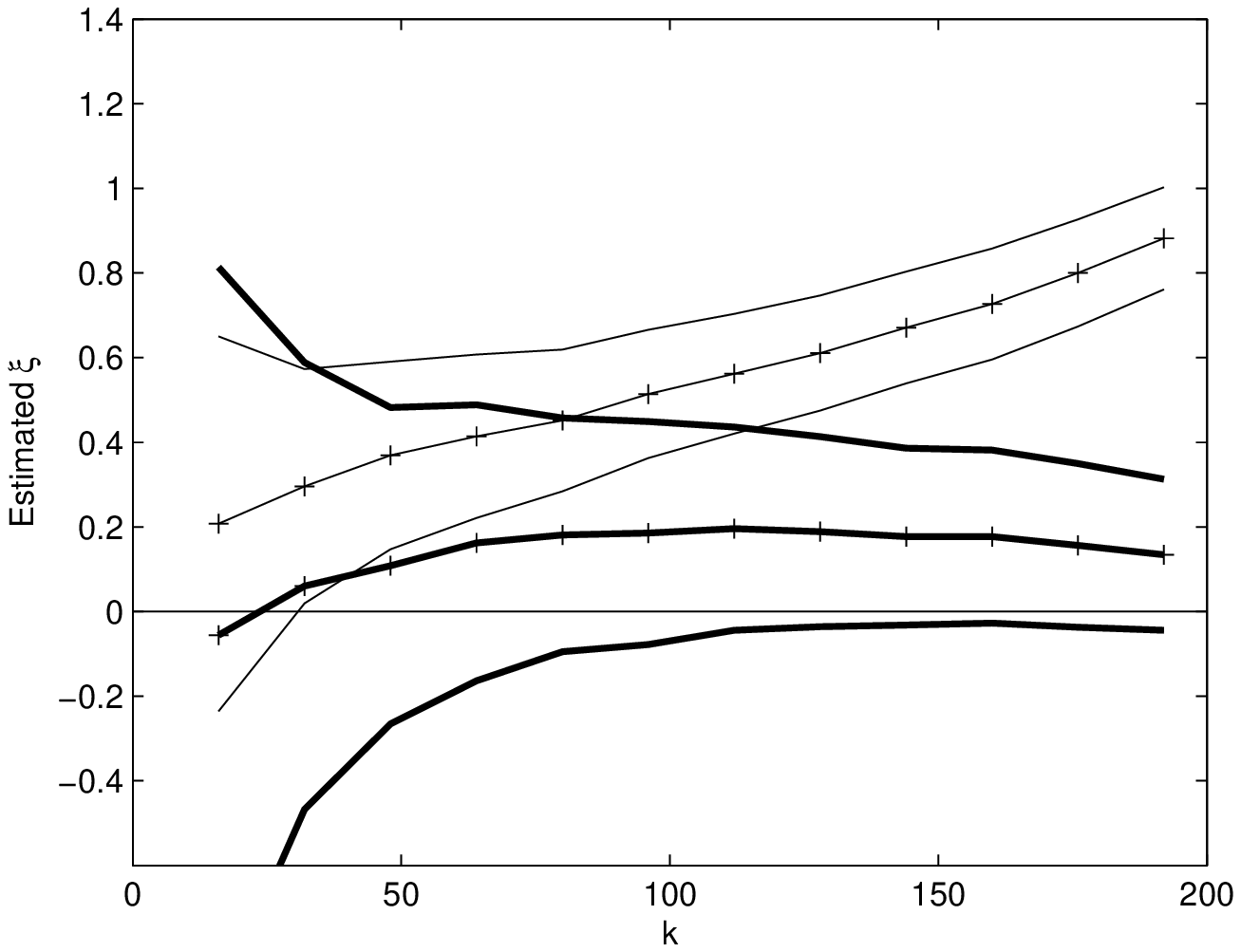}\\
\includegraphics[height = 60mm,keepaspectratio]{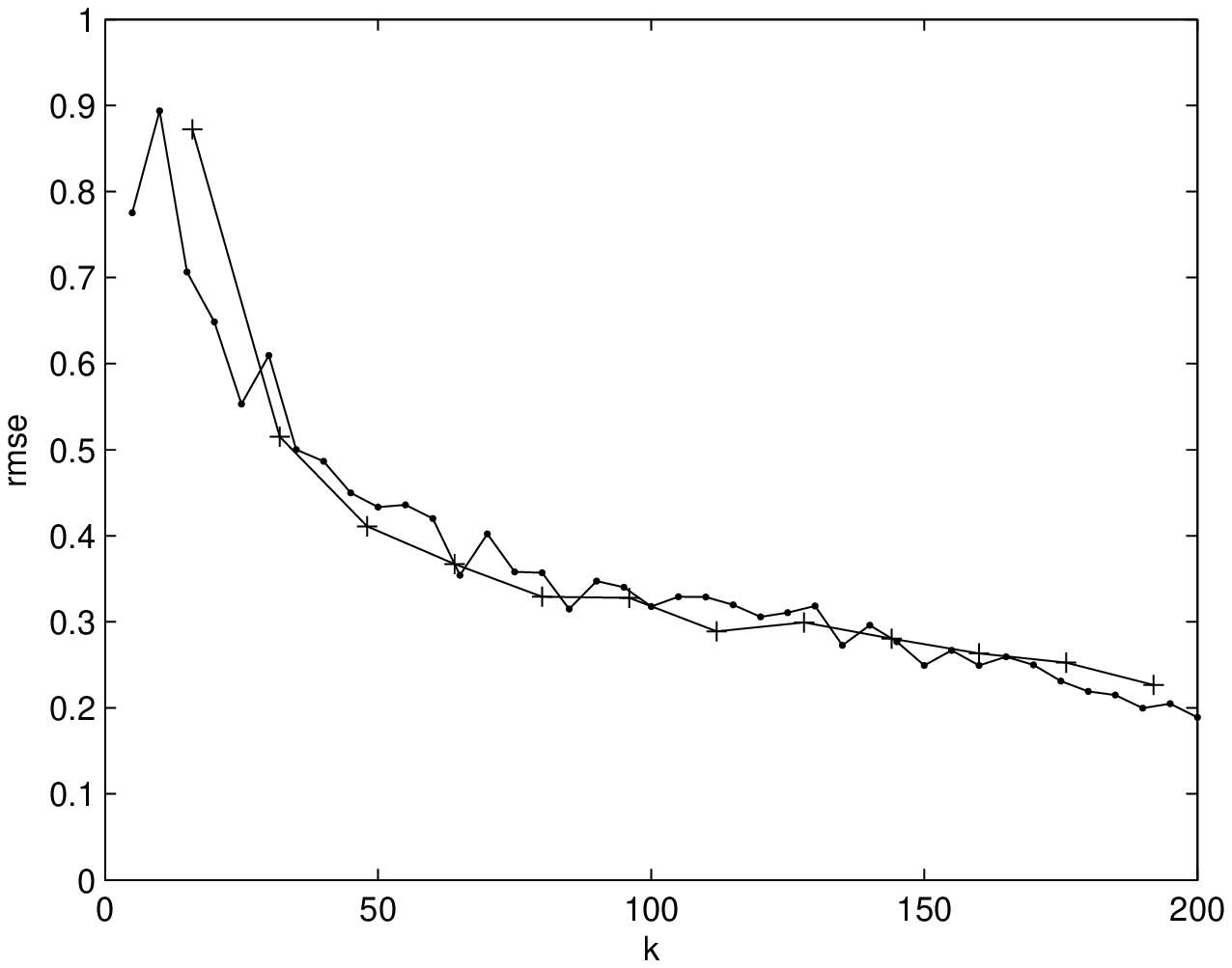}
  \caption{ Estimates (mean $\pm$ std in the top figures, RMSE in the lower figure) of the tail index of $\tau = 0.5$ Weibull samples of size $N=200$ using only the $k$ upper order statistics.
  The Segers estimates are top left, and the curve-fit estimates are top right.
  For the Segers estimates, the initial unconstrained estimates are shown by the thinner lines, with the thicker lines
  corresponding to the asymptotically unbiased constrained estimates.
  For the curve fit estimates, the thinner lines correspond simply fitting a curve to the first $k/2$ order statistics (normalised wrt the ($k/2$)th and $k$th), whilst the thicker lines make eight curve-fits at each $k$ (corresponding to $k'/k = [1:8]/8$), with the final estimate being a weighted least squares linear back-extrapolation to the hypothetical $k'=0$.
  The lower figure shows the RMSE for the constrained Segers (.) and the $k'$ curve fit (+) to be close.
  }\label{finalcomp}
\end{figure}

A simple procedure is outlined here for demonstration purposes only. For a sample of size $N$, for any $k \leq N$, we can pick the length-eight sequence $k' = \text{round}((m/8)k), m = 1, \ldots ,8$ at the eighth points of $k$. We can then construct eight curve-fit estimates, using only the first $k'/2$ order statistics normalised with respect to the $k'/2$th and $k'$th values (and using the plotting positions $G_i = (i-0.5)/k$ ). The eight estimates corresponding to the eight $k'$ can then be extrapolated back to the hypothetical $k'=0$ intercept by least squares fitting of a straight line (and, for this example, we choose to apply weights proportional to $(k')^2$ to give more weight to larger $k'$). The resulting estimates for the Weibull with $\tau = 0.5$ are shown in Figure~\ref{finalcomp}.

\cite{segers} proposed two versions of his tail index estimator. The unconstrained version was based heavily on the GPD, and was thus prone to substantial bias when applied to non-GPD data (as shown by Figure~\ref{segersweib}). The more sophisticated constrained version endeavoured to remove this bias by, loosely speaking, attempting to measure how non-GPD the data was and compensate accordingly. There is an obvious parallel with the $k$ and $k'$ curve-fit estimators here. The simple $k$-based curve-fit is constructed for the GPD, and the estimator that uses the sequence of $k'$ endeavours to compensate for non-GPD behaviour.

Figure~\ref{finalcomp} shows the corresponding performances of the unconstrained/constrained Segers and the $k/k'$ curve fit estimators to be strikingly similar. The lower part of the Figure shows the RMSE results for the enhanced estimators that have endeavoured to compensate for non-GPD behaviour. Considering the simplicity of the curve fit approach, the close correspondence is perhaps surprising.

The example was for illustration only. On some distributions (e.g. the normal) the constrained Segers gives lower RMSE than the $k'$ curve-fit, but on others (e.g. the (2,2) Beta distribution) the $k'$ curve-fit sequence is arguably the better. Although there is scope for further optimisation of the $k'$ estimator, the intention here was never to create a better tail index estimator, but to show how the curve-fit gives insight. The emphasis of the main body of the paper is on prediction: estimation is of lesser importance, especially since prediction is done using only the $k=20$ information, and at $k=20$ even the highly sophisticated constrained Segers estimator tends to give considerably worse RMSE than many simpler alternatives.

\end{document}